\documentclass[reqno]{amsart}

\usepackage{amssymb,amsmath,amsthm,xypic,amsfonts,latexsym,booktabs,qtree,tikz}

\theoremstyle{definition}
\newtheorem{definition}{Definition}

\newtheorem{remark}[definition]{Remark}

\theoremstyle{plain}

\newtheorem{theorem}[definition]{Theorem}

\newcommand{\ltree}[2]{\Tree[ .#1 [ .{\\[-7pt] #2} {\;} {\;} ] {\;} ]}
\newcommand{\rtree}[2]{\Tree[ .#1 {\;} [ .{\\[-7pt] #2} {\;} {\;} ] ]}

\newcommand{\lltree}[3]{\Tree[ .#1 [ .{\\[-7pt] #2} [ .{\\[-7pt] #3} {\;} {\;} ] {\;} ] {\;} ]}
\newcommand{\lrtree}[3]{\Tree[ .#1 [ .{\\[-7pt] #2} {\;} [ .{\\[-7pt] #3} {\;} {\;} ] ] {\;} ]}
\newcommand{\rltree}[3]{\Tree[ .#1 {\;} [ .{\\[-7pt] #2} [ .{\\[-7pt] #3} {\;} {\;} ] {\;} ] ]}
\newcommand{\rrtree}[3]{\Tree[ .#1 {\;} [ .{\\[-7pt] #2} {\;} [ .{\\[-7pt] #3} {\;} {\;} ] ] ]}

\newcommand{\corolla}[3]{\Tree[ .#1 [ .{\\[-7pt] #2} {\;} {\;} ] [ .{\\[-7pt] #3} {\;} {\;} ] ]}

\newcommand{\treeingamma}[1]{\Tree[ .{#1\\[-10pt]} {\;} {\;} ]}
\newcommand{\treeop}[2]{\; \begin{array}{c} \\[#2 pt] #1 \end{array} \;}
\newcommand{\ingamma}[2]{\; \begin{array}{c} \\[#2 pt] #1 \end{array} \;}

\allowdisplaybreaks

\begin{document}

\title[G-S bases for dendriform algebras and quadri-algebras]{Gr\"obner-Shirshov bases for the non-symmetric operads of dendriform algebras and quadri-algebras}

\author[S.~Madariaga]{Sara Madariaga}
\thanks{Department of Mathematics and Statistics, University of Saskatchewan, Canada. \\
\indent Email address: madariaga@math.usask.ca. \\
\indent The research of the author was supported by a postdoctoral fellowship from PIMS (Pacific Institute for the Mathematical Sciences)
and the Spanish Ministerio de Ciencia e Innovaci\'on (MTM2010-18370-C04-03).}

%

\begin{abstract}
    In this paper we use the operadic framework to find Gr\"obner-Shirshov bases
    for the free quadri-algebra.
    We perform computations using the representation of the nonsymmetric operad
    by planar rooted trees in a very intuitive way.
    Gr\"obner-Shirshov bases for the free dendriform algebra are also found with this
    technique, simplifying the work by Chen and Wang \cite{ChenWang10}.
\end{abstract}

\dedicatory{To Chema}



\maketitle


\section{Introduction}

Gr\"obner-Shirshov bases are a very powerful tool to study ideals of free algebras.
They provide an algorithmic way of computing normal forms of elements in the quotient ring
and help describing a linear basis for the quotient as well as finding its dimension.
Operads encode information about classes of algebras.

The goal of this paper is to compute Gr\"obner-Shirshov bases for the free dendriform algebra and
the free quadri-algebra using the theory of operads. The results for dendriform algebras improve
and simplify the work by Chen and Wang \cite{ChenWang10} in the following way.
Our approach via operads provides us with a more intuitive proof of their result and we also found a shorter
Gr\"obner-Shirshov basis by changing the order of operations.

Some computations were performed with the computer algebra system Maple 16 using the algorithms developed
and implemented by the author but the results on this paper do not depend on Maple computations.

We start defining the structures appearing and their main properties.
Throughout this paper, $F$ denotes a field of characteristic zero.

\subsection*{Operads}

Let $\mathrm{P}$-alg be a category of algebras (or, in other words, a type of algebras). An object of $\mathrm{P}$-alg is a vector space $A$ endowed
with $n$-ary generating operations $\mu_i \colon A^{\otimes n} \to A$ for (possibly for various $n$) satisfying some multilinear relations $r_j=0$
of the form
\[
    \sum_\phi \phi(a_1, \dots, a_n) = 0 \text{ for all } a_1, \dots, a_n \in A,
\]
where $\phi$ is a composition of the generating operations $\mu_i$. An element like $r = \sum \phi$ is called a relator.
(Think for instance about the classical definition of associative algebras, commutative algebras, Lie algebras, etc.)

If we consider $\mathcal{P}(n)$ the vector space of the $n$-ary operations defined on a $\mathrm{P}$-algebra, then for each $n$ we have a linear map
\begin{align*}
    \phi \colon \mathcal{P}(n) \otimes A^{\otimes n} & \longrightarrow A \\
    ( \mu; x_1, \dots, x_n ) & \mapsto \mu(x_1, \dots, x_n).
\end{align*}
The symmetric group $S_n$ acts on $A^{\otimes n}$ on the left and on $\mathcal{P}(n)$ on the right and these actions are compatible:
\[
    \mu^\sigma(x_1, \dots, x_n) = \mu(\sigma \cdot (x_1, \dots, x_n))
\]
We can compose operations in the following way. Given two operations $\mu \in \mathrm{P}(m)$ and $\nu \in \mathrm{P}(n)$ and and index $i$,
$1 \leq i \leq m$, we define the composition $ \mu \circ_i \nu \in \mathrm{P}(m+n-1) $ by
\[ \scriptsize
    \mu \circ_i \nu( x_1, \dots, x_{m+n-1} ) = \mu( x_1, \dots, x_{i-1}, \nu(x_i, \dots, x_{i+n-1}),x_{i+n}, \dots, x_{m+n-1} ).
\]
Compositions are associative:
\begin{align*}
    ( \lambda \circ_i \mu ) \circ_{j+m-1} \nu = ( \lambda \circ_j \nu ) \circ_i \mu & \quad \text{ if } 1 \leq i < j \leq l, \\
    ( \lambda \circ_i \mu ) \circ_{i-1+j} \nu = \lambda \circ_i ( \mu \circ_i \nu ) & \quad \text{ if } 1 \leq i \leq l, 1 \leq j \leq m.
\end{align*}
A linear operad $\mathcal{P}$ is given by a collection of vector spaces $\{ \mathcal{P}(n) \}_{n \geq 0}$ with
corresponding actions of the symmetric groups $S_n$ and compositions $\circ_i$ verifying the properties above.

We denote $\mathbb{S} = ( S_0, S_1, S_2, \dots, S_k, \dots )$ and consider $M = (M_0, M_1, \dots, M_k, \dots )$ the $\mathbb{S}$-module
which is, in arity $n$, the $S_n$-module spanned by the generating $n$-ary operations
(the $S_n$-module structure is determined by the symmetries of the operations).

A relator determines an operation in the operad $\mathcal{T}M$ generated by $M$. Let $R$ be the sub-$\mathbb{S}$-module
of $\mathcal{T}M$ spanned by all relators and let $(R)$ be the operadic ideal of $\mathcal{T}M$ generated by $R$.
Then we can construct the quotient operad $\mathcal{T}M /(R)$.

If $\mathcal{P}(V)$ is the free $\mathrm{P}$-algebra over $V$ and we set $\mathcal{P}(n)$ to be the the space of multilinear operations of arity $n$ 
regarded as an $S_n$-module then $\{ \mathcal{P}(n) \}_{n \geq 0}$ is an algebraic operad and coincides with the quotient operad $\mathcal{T}M /(R)$.
Moreover, a type of algebras whose relations are multilinear determine an operad and the category of algebras over this operad is equivalent to the category of algebras of the given type.

This generalization of the notion of types of algebras allows an operadic approach to many aspects of their study, as homology and Koszul duality.

In this paper we work only with nonsymmetric operads, that is, algebraic operads such that
\[
    \mathcal{P}(V) = \bigoplus_{n \geq 0} \mathcal{P}_n \otimes_F V^{\otimes n}.
\]
We recall here that $V$ is the space of elements and $\mathcal{P}_n$ is the space of $n$-ary operations defined over $V$.
With this notation, $\mathcal{P}(n)=\mathcal{P}_n \otimes_F FS_n$.
This relation allows us to convert results about operads into results about free algebras.

These operads encode types of algebras for which the generating operations do not satisfy any symmetry properties, and the variables
stay in the same order in each monomial of the defining relations. Dendriform algebras and quadri-algebras are examples of
algebras giving rise to nonsymmetric operads.

The combinatorial objects involved in the description of the free nonsymmetric operad are the planar rooted trees, which we will
use for our computations. The internal vertices of a tree monomial are labeled by the operations of the corresponding monomial.
In this paper we deal only with binary operads. The correspondence between monomials and tree monomials works as follows for
degrees 2 and 3 (higher degrees are similar). \vspace{-25pt}
\[
    \hspace{-25pt}
    \treeop{ \treeop{\mu_i(x,y) \; = }{4} \tiny \Tree[ .$\mu_i$ {\;} {\;} ] }{10} \normalsize
    \hspace{-10pt}
    \treeop{\mu_i(\mu_j(x,y),z) \; = }{15} \tiny \ltree{$\mu_i$}{$\mu_j$} \normalsize
    \;
    \treeop{\mu_i(x,\mu_j(y,z)) \; = }{15} \tiny \rtree{$\mu_i$}{$\mu_j$}
\]
We are especially interested in the following set of operations defined on trees:
\[
    \gamma_{n,i_1,\dots,i_n} \colon \mathcal{T}(X)_n \times \mathcal{T}(X)_{i_1} \times \cdots \times \mathcal{T}(X)_{i_n} \to \mathcal{T}(X)_{i_1+\dots+i_n},
\]
where $\mathcal{T}_k$ is the set of all linear combinations of tree monomials with $k$ leaves.
Roughly speaking, these operations take a tree monomial with $n$ leaves and attach the other $n$ tree monomials at its leaves.

For a comprehensive monograph on algebraic operads, see the book \cite{LodayVallette12} by Loday and Vallette.

\subsection*{Dendriform (di)algebras}

Dendriform algebras first appeared (with a different terminology) in the work of Loday and Ronco \cite{LodayRonco98} and
were named by Loday \cite{Loday01} in his study of HY homology theory. A dendriform
algebra is a vector space endowed with two binary operations $\prec$ and $\succ$ (here is the reason for sometimes calling
them dendriform dialgebras) related by the identities
\begin{align}\label{dendriform}
    ( a \succ b ) \prec c & = a \succ ( b \prec c ), \notag \\
    ( a \prec b ) \prec c & = a \prec (b \prec c ) + a \prec ( b \succ c ), \\
    a \succ ( b \succ c ) & = ( a \prec v ) \succ c + ( a \succ b ) \prec c. \notag
\end{align}
Note that the monomials involved in the relations keep the variables in the same order.
This concept dichotomizes the notion of associative algebra because the product $ x * y = x \prec y + x \succ y $
defined on a dendriform algebra is associative.

There exist several generalizations of this idea of decomposing operations into sums of coherent operations,
namely the trialgebra of Loday and Ronco \cite{LodayRonco04}, the quadri-algebra of Aguiar and Loday \cite{AguiarLoday04},
the ennea-algebra, the NS-algebra, the dendriform-Nijenhuis algebra and the octo-algebra of Leroux \cite{Leroux03,Leroux04a,Leroux04b}.
All these structures can be reconstructed from known objects by means of linear operators, see \cite{EbrahimiFardGuo08,Guo04}.

Dendriform algebras are also closely related to associative dialgebras in the framework of Koszul duality for algebraic operads.

\subsection*{Quadri-algebras}

Quadri-algebras were introduced by Aguiar and Loday \cite{AguiarLoday04} in their study of the space of endomorphisms of an
infinitesimal bialgebra. They recalled the existence of a pair of commuting Baxter operators determining two dendriform structures
on the space of endomorphisms and wanted to see how they are intertwined.

A quadri-algebra is a vector space $Q$ endowed with four binary operations $\searrow$, $\nearrow$, $\nwarrow$ and $\swarrow$
satisfying the relations

\begin{align}\label{quadrialgebra}
    ( x \nwarrow y ) \nwarrow z & = x \nwarrow ( y * z )  & \quad
    ( x \nearrow y ) \nwarrow z & = x \nearrow ( y \prec z ) \notag \\*
    ( x \wedge y ) \nearrow z & = x \nearrow ( y \succ z ) & \quad
    ( x \swarrow y ) \nwarrow z & = x \swarrow ( y \wedge z) \notag \\*
    ( x \searrow y ) \nwarrow z & = x \searrow ( y \nwarrow z )  & \quad
    ( x \vee y ) \nearrow z & = x \searrow ( y \nearrow z )  \\*
    ( x \prec y ) \swarrow z & = x \swarrow ( y \vee z )  & \quad
    ( x \succ y ) \swarrow z & = x \searrow ( y \swarrow z )  \notag \\*
    ( x * y ) \searrow z & = x \searrow ( y \searrow z ) &  \notag
\end{align}
where
\begin{alignat*}{2}
    x \succ y = x \nearrow y + x \searrow y & \qquad x \prec y = x \nwarrow y + x \swarrow y \\*
    x \vee y = x \searrow y + x \swarrow y & \qquad x \wedge y = x \nearrow y + x \nwarrow y
\end{alignat*}
and
\begin{align*}
    x * y  & \; = \; x \succ y + x \prec y \; = \; x \searrow y + x \nearrow y + x \nwarrow y + x \swarrow y.
\end{align*}

Note that $(Q,\prec, \succ)$ and $(Q,\wedge,\vee)$ are dendriform algebras and $(Q,*)$ is an associative algebra.
So a quadri-algebra may be seen as an associative algebra for which the multiplication can be decomposed into the
sum of four coherent operations.

Quadri-algebras can also be obtained from dendriform algebras and a Baxter operator
or from an associative algebra and a pair of commuting Baxter operators.
Basic examples of quadri-algebras are the algebra of permutations and the shuffle algebra.

Note that the four operations defining a quadri-algebra have no symmetries, and the relations (\ref{quadrialgebra})
involve only monomials where the variables $x,y,z$ stay in the same order.
This implies that the corresponding operad is nonsymmetric.

\subsection*{Gr\"obner-Shirshov bases}

The study of algebras by means of presentations by generators and relations raises the following key questions:
how to decide if a given element belongs to the ideal generated by the relations or how to find a linear basis of such algebras.
Shirshov \cite{Shirshov,Shirshov62} and Buchberger \cite{Buchberger87} gave independently a similar algorithmic solution to these problems
for Lie and commutative algebras respectively. Shirshov also proved that a similar algorithm is also valid for associative algebras.

The main idea of the algorithm is the following.
Let $S$ be a subset of the corresponding free algebra $F \langle X \rangle$ on the set $X$. Let $S'$ be the set resulting of adding
all non-trivial compositions (after reduction by $S$) of elements of $S$ to $S$ (note that $S \subseteq S'$).
Iterating this process one gets a sequence $ S \subseteq S' \subseteq S'' \subseteq \cdots \subseteq S^{(n)} \subseteq \cdots $.
The union $S^c$ of the sequence verifies that any composition of elements of $S^c$ is trivial relative to $S^c$.
It is what is now called a Gr\"obner-Shirshov basis.
Then from the Composition-Diamond Lemma it follows that the set all $S^c$-irreducible monomials is a linear basis of the
algebra generated by $X$ with defining relations $S$. Also, a given polynomial belongs to the ideal generated by $S$ if and only if
it reduces to zero with respect to $S$.

Intuitively, the Composition-Diamond Lemma says that if we want to compute the normal form of an element $f \in F \langle X \rangle$ with
respect to a subset $S \subset F \langle X \rangle$, at each step we could have many choices of reduction by different elements of $S$
but the final result does not depend on the sequence of reductions we perform.
The Composition-Diamond Lemma has been proved for a number of structures, namely associative algebras, semisimple Lie (super) algebras, irreducible
modules, Kac-Moody algebras, Coxeter groups, braid groups, conformal algebras, Loday's dialgebras, Leibniz algebras, multioperator algebras,
Rota-Baxter algebras, operads,~... In particular, for non-symmetric operads a proof can be found in \cite[Theorem~2.1]{DotsenkoVallette}.
This allows the construction of Gr\"obner-Shirshov bases for these structures.

All these compositions and reductions are performed with respect to a monomial order. In this paper we use the path-lexicographical order of
tree monomials. To each tree monomial we associate a sequence containing for each leaf (leaves ordered from left to right as they appear in
the monomial) the sequence of vertices needed to get from the root to it. For example, the path sequence associated to the tree monomial \\[-15pt]
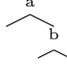
\[ \tiny{\rtree{a}{b}} \]
is $(a,ab,ab)$. To compare two path sequences $(w_1, \dots, w_{n_w})$ and $(w'_1, \dots, w'_{n_{w'}})$ we first compare their lengths
$n_w$ and $n_{w'}$ (which correspond to the degree of the monomials) and the longest sequence gives the greater tree monomial.
If $n_w = n_{w'}$ then we compare degree-lexicographically word by word of the sequence until we find a pair $(w_i,w'_i)$ which are different and
then the monomial corresponding to the greater sequence is greater. If not, we declare the tree monomials to be equal. For example: \vspace{-25pt}
\[
    \treeop{ \tiny \Tree[ .$\mu_i$ {\;} {\;} ] }{10} \normalsize
    \hspace{-20pt}
    \treeop{ \;\;\;\; < }{15} \tiny \ltree{$\mu_i$}{$\mu_j$} \normalsize
    \treeop{ < }{15} \tiny \rtree{$\mu_i$}{$\mu_j$}
\]

For a more detailed discussion on Gr\"obner-Shirshov bases, see \cite{BokutChen13}


\section{Gr\"obner-Shirshov bases for the dendriform operad}

Dendriform algebras are defined by the identities (\ref{dendriform}).
Rewriting these identities in terms of tree monomials gives the following polynomials equal to zero: \\[-25pt]
\begin{align*}
    & \treeop{(D1) \; =}{20} \tiny \ltree{$\prec$}{$\succ$} \treeop{-}{20} \rtree{$\succ$}{$\prec$} \\[-6pt]
    & \treeop{(D2) \; =}{20} \tiny \ltree{$\prec$}{$\prec$} \treeop{-}{20} \rtree{$\prec$}{$\prec$} \treeop{-}{20} \rtree{$\prec$}{$\succ$} \\[-6pt]
    & \treeop{(D3) \; =}{20} \tiny \rtree{$\succ$}{$\succ$} \treeop{-}{20} \ltree{$\succ$}{$\succ$} \treeop{-}{20} \ltree{$\succ$}{$\prec$} \treeop{.}{20}
\end{align*}

Note that since the operations defining dendriform algebras are binary, the generators of the corresponding operad are concentrated in arity 2.
Now we have to choose an ordering of monomials to determine the leading monomial of each identity and compute the compositions.
We use the path-lexicographical ordering of monomials together with a total order on the set of operations $\{ \prec, \succ \}$.
This choice is crucial: if we change the order of operations we obtain different results.
By setting $\prec$ precedes $\succ$ we obtain that the defining relations for dendriform algebras are a Gr\"obner-Shirshov basis of the dendriform operad
whereas if we assume $\succ$ precedes $\prec$ we recover the Gr\"obner-Shirshov basis found by Chen and Wang \cite{ChenWang10}, which has one more
element, resulting from a composition that cannot be reduced to zero.

\subsection{Ordering of operations `` $ \prec \textit{ precedes } \succ $ '' }

We assume that $\prec$ precedes $\succ$ in our set of operations, then identify the leading monomials of (D1), (D2), (D3)
with respect to this ordering and rewrite (if necessary) the identities so the leading monomial appears first and with coefficient 1: \\[-30pt]
\begin{align*}
    \treeop{S = \Big\{ (D1') \; = \; }{15} & \tiny \ltree{$\prec$}{$\succ$} \treeop{-}{20} \rtree{$\succ$}{$\prec$}\treeop{,}{20} \\[-20pt]
    \treeop{ (D2') \; = \; }{15} & \tiny \ltree{$\prec$}{$\prec$} \treeop{-}{20} \rtree{$\prec$}{$\prec$} \treeop{-}{20} \rtree{$\prec$}{$\succ$}\treeop{,}{20} \\[-20pt]
    \treeop{ (D3') \; = \; }{15} & \tiny \ltree{$\succ$}{$\succ$} \treeop{+}{20} \ltree{$\succ$}{$\prec$} \treeop{-}{20} \rtree{$\succ$}{$\succ$}
    \treeop{ \Bigg\} .}{15}
\end{align*}
\normalsize
Now we find the compositions of these polynomials.
\begin{itemize}
    \item[](D1')-(D1'): there are no compositions.
    \item[](D1')-(D2'): there is one composition:\\[-22pt]
        \begin{align*}
            \qquad\qquad
            \ingamma{\gamma_{2,3,1}\left(\right.}{4} \treeingamma{\\[-6pt]$\prec$} \ingamma{, \; (D1') \; , \; \vert \;\left.\right)}{4}
            \hspace{-15pt} \treeop{-}{4} \hspace{-15pt}
            \ingamma{\gamma_{3,2,1,1}\left(\right.\; (D2') \; ,}{4} \treeingamma{\\[-6pt]$\succ$} \ingamma{, \; \vert \; , \; \vert \;\left.\right)}{4}
        \end{align*}\\[-25pt]
    \item[](D1')-(D3'): there is one composition:\\[-22pt]
        \begin{align*}
            \qquad\qquad
            \ingamma{\gamma_{2,3,1}\left(\right.}{4} \treeingamma{\\[-6pt]$\prec$} \ingamma{, \; (D3') \; , \; \vert \;\left.\right)}{4}
            \hspace{-15pt} \treeop{-}{4} \hspace{-15pt}
            \ingamma{\gamma_{3,2,1,1}\left(\right.\; (D1') \; ,}{4} \treeingamma{\\[-6pt]$\succ$} \ingamma{, \; \vert \; , \; \vert \;\left.\right)}{4}
        \end{align*}\\[-25pt]
    \item[](D2')-(D2'): there is one composition:\\[-22pt]
        \begin{align*}
            \qquad\qquad
            \ingamma{\gamma_{2,3,1}\left(\right.}{4} \treeingamma{\\[-6pt]$\prec$} \ingamma{, \; (D2') \; , \; \vert \;\left.\right)}{4}
            \hspace{-15pt} \treeop{-}{4} \hspace{-15pt}
            \ingamma{\gamma_{3,2,1,1}\left(\right.\; (D2') \; ,}{4} \treeingamma{\\[-6pt]$\prec$} \ingamma{, \; \vert \; , \; \vert \;\left.\right)}{4}
        \end{align*}\\[-25pt]
    \item[](D2')-(D3'): there are no compositions.
    \item[](D3')-(D3'): there is one composition:\\[-22pt]
        \begin{align*}
            \qquad\qquad
            \ingamma{\gamma_{2,3,1}\left(\right.}{4} \treeingamma{\\[-6pt]$\succ$} \ingamma{, \; (D3') \; , \; \vert \;\left.\right)}{4}
            \hspace{-15pt} \treeop{-}{4} \hspace{-15pt}
            \ingamma{\gamma_{3,2,1,1}\left(\right.\; (D3') \; ,}{4} \treeingamma{\\[-6pt]$\succ$} \ingamma{, \; \vert \; , \; \vert \;\left.\right)}{4}
        \end{align*}\\[-25pt]
\end{itemize}
We now reduce the compositions with respect to $S$ to decide whether we have a Gr\"obner-Shirshov basis. We perform the reductions step by step using tree monomials
and indicate by a subscript the polynomial we use for the reductions in each step.\\[-22pt]
\begin{align*}
    & \ingamma{\gamma_{2,3,1}\left(\right.}{4} \treeingamma{\\[-6pt]$\prec$} \ingamma{, \; (D1') \; , \; \vert \;\left.\right)}{4}
    \hspace{-15pt} \treeop{-}{4} \hspace{-15pt}
    \ingamma{\gamma_{3,2,1,1} \left(\right.\; (D2') \; ,}{4} \treeingamma{\\[-6pt]$\succ$} \ingamma{, \; \vert \; , \; \vert \;\left.\right)}{4}
    \\[-15pt]
    & \quad \treeop{=}{25} \tiny
    \lltree{$\prec$}{$\prec$}{$\succ$}
    \treeop{-}{25}
    \lrtree{$\prec$}{$\succ$}{$\prec$}
    \treeop{-}{25}
    \lltree{$\prec$}{$\prec$}{$\succ$}
    \treeop{+}{25}
    \corolla{$\prec$}{$\succ$}{$\prec$}
    \treeop{+}{25}
    \corolla{$\prec$}{$\succ$}{$\succ$}
    \\[-15pt]
    & \quad \treeop{=}{25} \tiny
    \treeop{-}{25}
    \lrtree{$\prec$}{$\succ$}{$\prec$}_{(D1')}
    \treeop{+}{25}
    \corolla{$\prec$}{$\succ$}{$\prec$}_{(D1')}
    \treeop{+}{25}
    \corolla{$\prec$}{$\succ$}{$\succ$}_{(D1')}
    \\[-15pt]
    & \quad \treeop{\equiv}{25} \tiny
    \treeop{-}{25}
    \rltree{$\succ$}{$\prec$}{$\prec$}_{(D2')}
    \treeop{+}{25}
    \rrtree{$\succ$}{$\prec$}{$\prec$}
    \treeop{+}{25}
    \rrtree{$\succ$}{$\prec$}{$\succ$}
    \\[-15pt]
    & \quad \treeop{\equiv}{25} \tiny
    \treeop{-}{25}
    \rrtree{$\succ$}{$\prec$}{$\prec$}
    \treeop{-}{25}
    \rrtree{$\succ$}{$\prec$}{$\succ$}
    \treeop{+}{25}
    \rrtree{$\succ$}{$\prec$}{$\prec$}
    \treeop{+}{25}
    \rrtree{$\succ$}{$\prec$}{$\succ$}
    \normalsize \treeop{= 0}{25}
    \\
    & \ingamma{\gamma_{2,3,1}\left(\right.}{4} \treeingamma{\\[-6pt]$\prec$} \ingamma{, \; (D3') \; , \; \vert \;\left.\right)}{4}
    \hspace{-15pt} \treeop{-}{4} \hspace{-15pt}
    \ingamma{\gamma_{3,2,1,1}\left(\right.\; (D1') \; ,}{4} \treeingamma{\\[-6pt]$\succ$} \ingamma{, \; \vert \; , \; \vert \;\left.\right)}{4}
    \\[-15pt]
    & \quad \treeop{=}{25} \tiny
    \lltree{$\prec$}{$\succ$}{$\succ$}
    \treeop{-}{25}
    \lrtree{$\prec$}{$\succ$}{$\succ$}
    \treeop{+}{25}
    \lltree{$\prec$}{$\succ$}{$\prec$}
    \treeop{-}{25}
    \lltree{$\prec$}{$\succ$}{$\succ$}
    \treeop{+}{25}
    \corolla{$\succ$}{$\succ$}{$\prec$}
    \\[-15pt]
    & \quad \treeop{=}{25} \tiny
    \treeop{-}{25}
    \lrtree{$\prec$}{$\succ$}{$\succ$}_{(D1')}
    \treeop{+}{25}
    \lltree{$\prec$}{$\succ$}{$\prec$}_{(D1')}
    \treeop{+}{25}
    \corolla{$\succ$}{$\succ$}{$\prec$}_{(D3')}
    \\[-15pt]
    & \quad \treeop{\equiv}{25} \tiny
    \treeop{-}{25}
    \rltree{$\succ$}{$\prec$}{$\succ$}
    \treeop{+}{25}
    \corolla{$\succ$}{$\prec$}{$\prec$}
    \treeop{+}{25}
    \rrtree{$\succ$}{$\succ$}{$\prec$}
    \treeop{-}{25}
    \corolla{$\succ$}{$\prec$}{$\prec$}
    \\[-15pt]
    & \quad \treeop{=}{25} \tiny
    \treeop{-}{25}
    \rltree{$\succ$}{$\prec$}{$\succ$}_{(D1')}
    \treeop{+}{25}
    \rrtree{$\succ$}{$\succ$}{$\prec$}
    \treeop{\equiv}{25}
    \treeop{-}{25}
    \rrtree{$\succ$}{$\succ$}{$\prec$}
    \treeop{+}{25}
    \rrtree{$\succ$}{$\succ$}{$\prec$}
    \normalsize \treeop{= 0}{25}
    \\
    & \ingamma{\gamma_{2,3,1}\left(\right.}{4} \treeingamma{\\[-6pt]$\prec$} \ingamma{, \; (D2') \; , \; \vert \;\left.\right)}{4}
    \hspace{-15pt} \treeop{-}{4} \hspace{-15pt}
    \ingamma{\gamma_{3,2,1,1}\left(\right.\; (D2') \; ,}{4} \treeingamma{\\[-6pt]$\prec$} \ingamma{, \; \vert \; , \; \vert \;\left.\right)}{4}
    \\[-15pt]
    & \quad \treeop{=}{25} \tiny
    \lltree{$\prec$}{$\prec$}{$\prec$}
    \treeop{-}{25}
    \lrtree{$\prec$}{$\prec$}{$\prec$}
    \treeop{-}{25}
    \lrtree{$\prec$}{$\prec$}{$\succ$}
    \treeop{-}{25}
    \lltree{$\prec$}{$\prec$}{$\prec$}
    \treeop{+}{25}
    \corolla{$\prec$}{$\prec$}{$\prec$}
    \treeop{+}{25}
    \corolla{$\prec$}{$\prec$}{$\succ$}
    \\[-15pt]
    & \quad \treeop{=}{25} \tiny
    \treeop{-}{25}
    \lrtree{$\prec$}{$\prec$}{$\prec$}_{(D2')}
    \treeop{-}{25}
    \lrtree{$\prec$}{$\prec$}{$\succ$}_{(D2')}
    \treeop{+}{25}
    \corolla{$\prec$}{$\prec$}{$\prec$}_{(D2')}
    \treeop{+}{25}
    \corolla{$\prec$}{$\prec$}{$\succ$}_{(D2')}
    \\[-15pt]
    & \quad \treeop{\equiv}{25} \tiny
    \treeop{-}{25}
    \rltree{$\prec$}{$\prec$}{$\prec$}_{(D2')}
    \treeop{-}{25}
    \rltree{$\prec$}{$\succ$}{$\prec$}
    \treeop{-}{25}
    \rltree{$\prec$}{$\prec$}{$\succ$}_{(D1')}
    \treeop{-}{25}
    \rltree{$\prec$}{$\succ$}{$\succ$}_{(D3')}
    \\[-15pt]
    & \qquad \qquad \tiny
    \treeop{+}{25}
    \rrtree{$\prec$}{$\prec$}{$\prec$}
    \treeop{+}{25}
    \rrtree{$\prec$}{$\succ$}{$\prec$}
    \treeop{+}{25}
    \rrtree{$\prec$}{$\prec$}{$\succ$}
    \treeop{+}{25}
    \rrtree{$\prec$}{$\succ$}{$\succ$}
    \\[-15pt]
    & \quad \treeop{\equiv}{25} \tiny
    \treeop{-}{25}
    \rrtree{$\prec$}{$\prec$}{$\prec$}
    \treeop{-}{25}
    \rrtree{$\prec$}{$\prec$}{$\succ$}
    \treeop{-}{25}
    \rltree{$\prec$}{$\succ$}{$\prec$}
    \treeop{-}{25}
    \rrtree{$\prec$}{$\succ$}{$\prec$}
    \treeop{-}{25}
    \rrtree{$\prec$}{$\succ$}{$\succ$}
    \\[-15pt]
    & \qquad \qquad \tiny
    \treeop{+}{25}
    \rltree{$\prec$}{$\succ$}{$\prec$}
    \treeop{+}{25}
    \rrtree{$\prec$}{$\prec$}{$\prec$}
    \treeop{+}{25}
    \rrtree{$\prec$}{$\succ$}{$\prec$}
    \treeop{+}{25}
    \rrtree{$\prec$}{$\prec$}{$\succ$}
    \treeop{+}{25}
    \rrtree{$\prec$}{$\succ$}{$\succ$}
    \normalsize \treeop{= 0}{25}
    \\
    & \ingamma{\gamma_{2,3,1}\left(\right.}{4} \treeingamma{\\[-6pt]$\succ$} \ingamma{, \; (D3') \; , \; \vert \;\left.\right)}{4}
    \hspace{-15pt} \treeop{-}{4} \hspace{-15pt}
    \ingamma{\gamma_{3,2,1,1}\left(\right.\; (D3') \; ,}{4} \treeingamma{\\[-6pt]$\succ$} \ingamma{, \; \vert \; , \; \vert \;\left.\right)}{4}
    \\[-15pt]
    & \quad \treeop{=}{25} \tiny
    \lltree{$\succ$}{$\succ$}{$\succ$}
    \treeop{-}{25}
    \lrtree{$\succ$}{$\succ$}{$\succ$}
    \treeop{+}{25}
    \lltree{$\succ$}{$\succ$}{$\prec$}
    \treeop{-}{25}
    \lltree{$\succ$}{$\succ$}{$\succ$}
    \treeop{+}{25}
    \corolla{$\succ$}{$\succ$}{$\succ$}
    \treeop{-}{25}
    \lltree{$\succ$}{$\prec$}{$\succ$}
    \\[-15pt]
    & \quad \treeop{=}{25} \tiny
    \treeop{-}{25}
    \lrtree{$\succ$}{$\succ$}{$\succ$}_{(D3')}
    \treeop{+}{25}
    \lltree{$\succ$}{$\succ$}{$\prec$}_{(D3')}
    \treeop{+}{25}
    \corolla{$\succ$}{$\succ$}{$\succ$}_{(D3')}
    \treeop{-}{25}
    \lltree{$\succ$}{$\prec$}{$\succ$}_{(D1')}
    \\[-15pt]
    & \quad \treeop{\equiv}{25} \tiny
    \treeop{-}{25}
    \rltree{$\succ$}{$\succ$}{$\succ$}
    \treeop{+}{25}
    \lrtree{$\succ$}{$\prec$}{$\succ$}
    \treeop{+}{25}
    \corolla{$\succ$}{$\prec$}{$\succ$}
    \treeop{-}{25}
    \lltree{$\succ$}{$\prec$}{$\prec$}
    \treeop{+}{25}
    \rrtree{$\succ$}{$\succ$}{$\succ$}
    \\[-15pt]
    & \qquad \qquad \tiny
    \treeop{-}{25}
    \corolla{$\succ$}{$\prec$}{$\succ$}
    \treeop{-}{25}
    \lrtree{$\succ$}{$\succ$}{$\prec$}
    \\[-15pt]
    & \quad \treeop{=}{25} \tiny
    \treeop{-}{25}
    \rltree{$\succ$}{$\succ$}{$\succ$}_{(D3')}
    \treeop{+}{25}
    \lrtree{$\succ$}{$\prec$}{$\succ$}
    \treeop{-}{25}
    \lltree{$\succ$}{$\prec$}{$\prec$}_{(D2')}
    \treeop{+}{25}
    \rrtree{$\succ$}{$\succ$}{$\succ$}
    \treeop{-}{25}
    \lrtree{$\succ$}{$\succ$}{$\prec$}_{(D3')}
    \\[-15pt]
    & \quad \treeop{\equiv}{25} \tiny
    \treeop{-}{25}
    \rrtree{$\succ$}{$\succ$}{$\succ$}
    \treeop{+}{25}
    \rltree{$\succ$}{$\succ$}{$\prec$}
    \treeop{+}{25}
    \lrtree{$\succ$}{$\prec$}{$\succ$}
    \treeop{-}{25}
    \lrtree{$\succ$}{$\prec$}{$\prec$}
    \treeop{-}{25}
    \lrtree{$\succ$}{$\prec$}{$\succ$}
    \\[-15pt]
    & \qquad \qquad \tiny
    \treeop{+}{25}
    \rrtree{$\succ$}{$\succ$}{$\succ$}
    \treeop{-}{25}
    \rltree{$\succ$}{$\succ$}{$\prec$}
    \treeop{+}{25}
    \lrtree{$\succ$}{$\prec$}{$\prec$}
    \normalsize \treeop{= 0}{25}
\end{align*}
Since all compositions reduce to zero, we can state the following result.

\begin{theorem}
    The set of defining identities for dendriform algebras is a self-reduced Gr\"obner-Shirshov basis for the non-symmetric
    dendriform operad with respect to the path-lexicographical ordering of monomials and the order of operations
    $ \prec \; < \; \succ $.
\end{theorem}

Since the dendriform operad is quadratic and its Gr\"obner-Shirshov basis is also quadratic,
from \cite[Prop.~3.10]{Hoffbeck10} and \cite[Corollary~3]{DotsenkoKhoroshkin10}, we immediately get
that the non-symmetric dendriform operad is Koszul. Thus we found an alternative proof of this fact.

We can describe the normal tree monomials with respect to this Gr\"obner-Shirshov basis.
For normal tree monomials, any growth to the right is allowed, but the only permitted growths to the left are of type \\[-22pt]
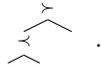
\[ {\tiny \ltree{$\succ$}{$\prec$} } \treeop{.}{15} \]

\begin{remark}
    Note that the polynomials in $S$ are simpler that the polynomials found by Chen and Wang \cite[Theorem 3.2]{ChenWang10}.
    Their approach to compute a basis of the free dendriform algebra as by considering it as a quotient of the free L-algebra.
    Non-symmetric operads possess a much richer structure (all compositions $\gamma$) than L-algebras (they have only two products),
    and it translates into a smaller set of restrictions on normal forms.
\end{remark}

Since the dimension of the multilinear subspaces of the free dendriform algebra is given by the Catalan number \cite{Loday01},
we have another description of this famous sequence of numbers in terms of rooted labeled planar binary trees with the previous
growth restrictions. This description is not part of Stanley's extensive list of different characterizations of the Catalan numbers
\cite{Stanley12}. A direct proof of this fact is an open problem.

\subsection{Ordering of operations `` $ \succ \textit{ precedes } \prec $ '' }

We assume that $\succ$ precedes $\prec$ in our set of operations, then identify the leading monomials of (D1), (D2), (D3)
with respect to this ordering and rewrite (if necessary) the identities so the leading monomial appears first and with coefficient 1: \\[-22pt]
\begin{align*}
    \treeop{S = \Big\{ (D1') \; = \; }{15} & \tiny \ltree{$\prec$}{$\succ$} \treeop{-}{20} \rtree{$\succ$}{$\prec$}\treeop{,}{20} \\[-20pt]
    \treeop{ (D2') \; = \; }{15} & \tiny \ltree{$\prec$}{$\prec$} \treeop{-}{20} \rtree{$\prec$}{$\prec$} \treeop{-}{20} \rtree{$\prec$}{$\succ$}\treeop{,}{20} \\[-20pt]
    \treeop{ (D3') \; = \; }{15} & \tiny \ltree{$\succ$}{$\prec$} \treeop{+}{20} \ltree{$\succ$}{$\succ$} \treeop{-}{20} \rtree{$\succ$}{$\succ$}
    \treeop{\Bigg\}.}{15}
\end{align*}
\normalsize
Now we find the compositions of these polynomials.
\begin{itemize}
    \item[](D1')-(D1'): there are no compositions.
    \item[](D1')-(D2'): there is one composition:\\[-22pt]
        \begin{align*}
            \qquad\qquad
            \ingamma{\gamma_{2,3,1}\left(\right.}{4} \treeingamma{\\[-6pt]$\prec$} \ingamma{, \; (D1') \; , \; \vert \;\left.\right)}{4}
            \hspace{-15pt} \treeop{-}{4} \hspace{-15pt}
            \ingamma{\gamma_{3,2,1,1}\left(\right.\; (D2') \; ,}{4} \treeingamma{\\[-6pt]$\succ$} \ingamma{, \; \vert \; , \; \vert \;\left.\right)}{4}
        \end{align*}\\[-25pt]
    \item[](D1')-(D3'): there are two compositions:\\[-22pt]
        \begin{align*}
            \qquad\qquad
            \ingamma{\gamma_{2,3,1}\left(\right.}{4} \treeingamma{\\[-6pt]$\succ$} \ingamma{, \; (D1') \; , \; \vert \;\left.\right)}{4}
            \hspace{-15pt} \treeop{-}{4} \hspace{-15pt}
            \ingamma{\gamma_{3,2,1,1}\left(\right.\; (D3') \; ,}{4} \treeingamma{\\[-6pt]$\succ$} \ingamma{, \; \vert \; , \; \vert \;\left.\right)}{4} \\[-20pt]
            \qquad\qquad
            \ingamma{\gamma_{2,3,1}\left(\right.}{4} \treeingamma{\\[-6pt]$\prec$} \ingamma{, \; (D3') \; , \; \vert \;\left.\right)}{4}
            \hspace{-15pt} \treeop{-}{4} \hspace{-15pt}
            \ingamma{\gamma_{3,2,1,1}\left(\right.\; (D1') \; ,}{4} \treeingamma{\\[-6pt]$\prec$} \ingamma{, \; \vert \; , \; \vert \;\left.\right)}{4}
        \end{align*}\\[-25pt]
    \item[](D2')-(D2'): there is one composition:\\[-22pt]
        \begin{align*}
            \qquad\qquad
            \ingamma{\gamma_{2,3,1}\left(\right.}{4} \treeingamma{\\[-6pt]$\prec$} \ingamma{, \; (D2') \; , \; \vert \;\left.\right)}{4}
            \hspace{-15pt} \treeop{-}{4} \hspace{-15pt}
            \ingamma{\gamma_{3,2,1,1}\left(\right.\; (D2') \; ,}{4} \treeingamma{\\[-6pt]$\prec$} \ingamma{, \; \vert \; , \; \vert \;\left.\right)}{4}
        \end{align*}\\[-25pt]
    \item[](D2')-(D3'): there is one composition:\\[-22pt]
        \begin{align*}
            \qquad\qquad
            \ingamma{\gamma_{2,3,1}\left(\right.}{4} \treeingamma{\\[-6pt]$\succ$} \ingamma{, \; (D2') \; , \; \vert \;\left.\right)}{4}
            \hspace{-15pt} \treeop{-}{4} \hspace{-15pt}
            \ingamma{\gamma_{3,2,1,1}\left(\right.\; (D3') \; ,}{4} \treeingamma{\\[-6pt]$\prec$} \ingamma{, \; \vert \; , \; \vert \;\left.\right)}{4}
        \end{align*}\\[-25pt]
    \item[](D3')-(D3'): there are no compositions.
\end{itemize}
If we reduce the compositions with respect to $S$ as before, we see that the last one does not reduce to zero,
so $S$ is not a Gr\"obner-Shirshov basis for the non-symmetric dendriform operad.
We now add the nonzero reduced compositions to $S$. We sort its terms and find its leading monomial: \\[-22pt]
\[
    \treeop{S' = S \cup \Big\{ \; (D4') = \!}{20}
    \tiny \lltree{$\succ$}{$\succ$}{$\succ$}
    \treeop{-}{20} \corolla{$\succ$}{$\succ$}{$\succ$}
    \treeop{+}{20} \lrtree{$\succ$}{$\succ$}{$\prec$}
    \treeop{ \Bigg\}.}{20}
\]
Now we need to perform the compositions of the new element with all the elements of $S'$.
\begin{itemize}
    \item[](D1')-(D4'): \\[-22pt]
        \begin{align*}
            \qquad\qquad
            \ingamma{\gamma_{2,4,1}\left(\right.}{4}
            \treeingamma{\\[-6pt]$\prec$}
            \ingamma{, \; (D4') \; , \; \vert \;\left.\right)}{4}
            \hspace{-15pt} \treeop{-}{4} \hspace{-15pt}
            \ingamma{\gamma_{3,3,1,1}\left(\right.\; (D1') \; ,}{4}
            {\tiny \Tree[ .{\\[-10pt]$\succ$} [ .{\\[-8pt]$\succ$} {\;} {\;} ] {\;} ]}
            \ingamma{, \; \vert \; , \; \vert \;\left.\right)}{4}
        \end{align*}\\[-25pt]
    \item[](D3')-(D4'): \\[-22pt]
        \begin{align*}
            \qquad\qquad
            \ingamma{\gamma_{3,3,1,1}\left(\right.}{4}
            {\tiny \Tree[ .{\\[-10pt]$\succ$} [ .{\\[-8pt]$\succ$} {\;} {\;} ] {\;} ]}
            \ingamma{, \; (D3') \; , \; \vert \; , \; \vert \;\left.\right)}{4}
            \hspace{-15pt} \treeop{-}{4} \hspace{-15pt}
            \ingamma{\gamma_{4,2,1,1,1}\left(\right.\; (D4') \; ,}{4}
            \treeingamma{\\[-6pt]$\prec$}
            \ingamma{, \; \vert \; , \; \vert \; , \; \vert \; \left.\right)}{4}
        \end{align*}\\[-25pt]
    \item[](D4')-(D4'): \\[-22pt]
        \begin{align*}
            \ingamma{\gamma_{2,4,1}\left(\right.}{4}
            \treeingamma{\\[-6pt]$\succ$}
            \ingamma{, \; (D4') \; , \; \vert \;\left.\right)}{4}
            &
            \hspace{-15pt} \treeop{-}{4} \hspace{-15pt}
            \ingamma{\gamma_{4,2,1,1,1}\left(\right.\; (D4') \; ,}{4}
            \treeingamma{\\[-6pt]$\succ$}
            \ingamma{, \; \vert \; , \; \vert \; , \; \vert \; \left.\right)}{4}
             \\[-20pt]
            \ingamma{\gamma_{3,4,1,1}\left(\right.}{4}
            {\tiny \Tree[ .{\\[-10pt]$\succ$} [ .{\\[-8pt]$\succ$} {\;} {\;} ] {\;} ]}
            \ingamma{, \; (D4') \; , \; \vert \; , \; \vert \; \left.\right)}{4}
            &
            \hspace{-15pt} \treeop{-}{4} \hspace{-15pt}
            \ingamma{\gamma_{4,3,1,1,1}\left(\right.\; (D4') \; ,}{4}
            {\tiny \Tree[ .{\\[-10pt]$\succ$} [ .{\\[-8pt]$\succ$} {\;} {\;} ] {\;} ]}
            \ingamma{, \; \vert \; , \; \vert \; , \; \vert \; \left.\right)}{4}
        \end{align*}\\[-25pt]
\end{itemize}
And all these compositions reduce to zero with respect to $S'$, so $S'$ is a Gr\"obner-Shirshov basis for
the non-symmetric dendriform operad.
This is exactly the Gr\"obner-Shirshov basis obtained by Chen and Wang in \cite{ChenWang10}.
Note that in this case the Gr\"obner-Shirshov basis is not quadratic because (D4') is cubic (it involves 3 operations).

We can describe the normal tree monomials with respect to this Gr\"obner-Shirshov basis.
For normal tree monomials, any growth to the right is allowed, but the only permitted growths to the left are of type \\[-22pt]
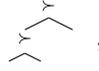
\[ {\tiny \ltree{$\succ$}{$\succ$} } \treeop{,}{15} \]
also excluding the possibility of having \\[-22pt]
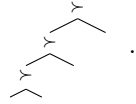
\[ {\tiny \lltree{$\succ$}{$\succ$}{$\succ$} } \treeop{.}{20} \]

In this way we get another description of the Catalan numbers in terms of rooted planar binary trees with these growth restrictions.


\section{Gr\"obner-Shirshov bases for the quadri-algebra operad}

Quadri-algebras are defined by the identities (\ref{quadrialgebra}). Since all the operations in these identities are binary,
the generators of the corresponding operad are concentrated in arity 2. 
Now we have to choose an ordering of monomials to determine the leading monomial of each identity and compute the compositions.
We use the path-lexicographical ordering of monomials together with a total order on the set of operations $\{ a,b,c,d \}$.
As seen in the previous section, this choice is important: if we change the order of operations we may obtain different results.
For quadri-algebras there are two orders of variables ($c<b<d<a$ and $c<d<b<a$) for which the defining identities form a
Gr\"obner-Shirshov basis of the corresponding non-symmetric operad, so it is Koszul. For some orders we obtain a Gr\"obner-Shirshov
basis after one iteration and for some others the number of non-zero reductions gets so big that keep computing is not worth the effort
(we already found a nice Gr\"obner-Shirshov basis for this non-symmetric operad).

Note that the set of tree monomials is symmetric in $b$ and $d$: if we interchange them then we get the same set.
This reduces the number of permutations of the set of operations we have to consider.
Table \ref{tablepermutations} displays the number of compositions and nonzero reductions for the different orderings of operations.
Whenever we find a zero in a \textit{red} column we have a Gr\"obner-Shirshov basis.

\begin{table}[h]
\caption{Gr\"obner bases for different orders of operations}
\label{tablepermutations}
\begin{tabular}{|c|cc|cc|cc|}
\toprule
    order   & \multicolumn{2}{|c|}{iteration 1} & \multicolumn{2}{|c|}{iteration 2} & \multicolumn{2}{|c|}{iteration 3} \\
               & \hspace{5pt} comp \hspace{5pt} & \hspace{5pt} red \hspace{5pt}
               & \hspace{5pt} comp \hspace{5pt} & \hspace{5pt} red \hspace{5pt}
               & \hspace{5pt} comp \hspace{5pt} & \hspace{5pt} red \hspace{5pt}     \\
\midrule
   \hspace{5pt} $a<b<c<d$ \hspace{5pt}   &        21       &       5       &        38       &      12       &        213      &     ??? \\
   \hspace{5pt} $a<b<d<c$ \hspace{5pt}   &        25       &      10       &        82       &      41       &       1411      &     ??? \\
   \hspace{5pt} $a<c<b<d$ \hspace{5pt}   &        19       &       3       &        20       &       7       &         80      &     ??? \\
   \hspace{5pt} $a<c<d<b$ \hspace{5pt}   &        19       &       3       &        20       &       7       &         80      &     ??? \\
   \hspace{5pt} $a<d<b<c$ \hspace{5pt}   &        25       &      10       &        82       &      41       &       1411      &     ??? \\
   \hspace{5pt} $a<d<c<b$ \hspace{5pt}   &        21       &       5       &        38       &      12       &        213      &     ??? \\
   \hspace{5pt} $b<a<c<d$ \hspace{5pt}   &        20       &       4       &        21       &       0       &         --      &      -- \\
   \hspace{5pt} $b<a<d<c$ \hspace{5pt}   &        24       &       9       &        66       &      23       &        970      &     ??? \\
   \hspace{5pt} $b<c<a<d$ \hspace{5pt}   &        20       &       4       &        21       &       0       &         --      &      -- \\
   \hspace{5pt} $b<c<d<a$ \hspace{5pt}   &        18       &       2       &        10       &       0       &         --      &      -- \\
   \hspace{5pt} $b<d<a<c$ \hspace{5pt}   &        23       &       8       &        62       &       0       &         --      &      -- \\
   \hspace{5pt} $b<d<c<a$ \hspace{5pt}   &        20       &       4       &        14       &       0       &         --      &      -- \\
   \hspace{5pt} $c<a<b<d$ \hspace{5pt}   &        19       &       3       &        19       &       4       &         49      &     ??? \\
   \hspace{5pt} $c<a<d<b$ \hspace{5pt}   &        19       &       3       &        19       &       4       &         49      &     ??? \\
   \hspace{5pt} $c<b<a<d$ \hspace{5pt}   &        18       &       2       &        10       &       0       &         --      &      -- \\
   \hspace{5pt} $c<b<d<a$ \hspace{5pt}   &        16       &       0       &        --       &      --       &         --      &      -- \\
   \hspace{5pt} $c<d<a<b$ \hspace{5pt}   &        18       &       2       &        10       &       0       &         --      &      -- \\
   \hspace{5pt} $c<d<b<a$ \hspace{5pt}   &        16       &       0       &        --       &      --       &         --      &      -- \\
   \hspace{5pt} $d<a<b<c$ \hspace{5pt}   &        24       &       9       &        66       &      23       &        970      &     ??? \\
   \hspace{5pt} $d<a<c<b$ \hspace{5pt}   &        20       &       4       &        21       &       0       &         --      &      -- \\
   \hspace{5pt} $d<b<a<c$ \hspace{5pt}   &        23       &       8       &        62       &       0       &         --      &      -- \\
   \hspace{5pt} $d<b<c<a$ \hspace{5pt}   &        20       &       4       &        24       &       0       &         --      &      -- \\
   \hspace{5pt} $d<c<a<b$ \hspace{5pt}   &        20       &       4       &        12       &       0       &         --      &      -- \\
   \hspace{5pt} $d<c<b<a$ \hspace{5pt}   &        18       &       2       &        10       &       0       &         --      &      -- \\
\bottomrule
\end{tabular}
\end{table}
We explicitly compute a Gr\"obner-Shirshov basis for the non-symmetric quadri-algebra operad for the order of operations
$c<b<d<a$ using tree monomials.
We first identify the leading monomials of (Q1)--(Q9) with respect to this ordering and rewrite
(if necessary) the identities so the leading monomial appears first and with coefficient 1: \\[-30pt]
\begin{align*}
    \treeop{S = \Big\{ (Q1') \; = \; }{15} &
    \tiny \ltree{c}{c} \treeop{-}{20} \rtree{c}{a} \treeop{-}{20} \rtree{c}{b} \treeop{-}{20} \rtree{c}{c} \treeop{-}{20} \rtree{c}{d}
    \treeop{,}{20} \\[-20pt]
    \treeop{ (Q2') \; = \; }{15} &
    \tiny \ltree{c}{d} \treeop{-}{20} \rtree{d}{b} \treeop{-}{20} \rtree{d}{c}
    \treeop{,}{20} \\[-20pt]
    \treeop{ (Q3') \; = \; }{15} &
    \tiny \ltree{d}{d} \treeop{+}{20} \ltree{c}{d} \treeop{-}{20} \rtree{d}{a} \treeop{-}{20} \rtree{d}{d}
    \treeop{,}{20} \\[-20pt]
    \treeop{ (Q4') \; = \; }{15} &
    \tiny \ltree{c}{b} \treeop{-}{20} \rtree{b}{c} \treeop{-}{20} \rtree{b}{d}
    \treeop{,}{20} \\[-20pt]
    \treeop{ (Q5') \; = \; }{15} &
    \tiny \ltree{c}{a} \treeop{-}{20} \rtree{a}{c}
    \treeop{,}{20} \\[-20pt]
    \treeop{ (Q6') \; = \; }{15} &
    \tiny \ltree{d}{a} \treeop{+}{20} \ltree{d}{b} \treeop{-}{20} \rtree{a}{d}
    \treeop{,}{20} \\[-20pt]
    \treeop{ (Q7') \; = \; }{15} &
    \tiny \ltree{b}{b} \treeop{+}{20} \ltree{b}{c} \treeop{-}{20} \rtree{b}{b} \treeop{-}{20} \rtree{b}{a}
    \treeop{,}{20} \\[-20pt]
    \treeop{ (Q8') \; = \; }{15} &
    \tiny \ltree{b}{a} \treeop{+}{20} \ltree{b}{d} \treeop{-}{20} \rtree{a}{b}
    \treeop{,}{20} \\[-20pt]
    \treeop{ (Q9') \; = \; }{15} &
    \tiny \ltree{a}{a} \treeop{+}{20} \ltree{a}{b} \treeop{+}{20} \ltree{a}{c} \treeop{+}{20} \ltree{a}{d} \treeop{-}{20} \rtree{a}{a}
    \treeop{\Bigg\}.}{15}
\end{align*}
\normalsize

Now we find the compositions of these polynomials. \\[-22pt]
\begin{align*}
    \ingamma{(Q1')-(Q1'):}{4} & \ingamma{\gamma_{2,3,1}\left(\right.}{4}
                                \hspace{-5pt}\treeingamma{\\[-6pt]c}\hspace{-10pt}
                                \ingamma{, \; (Q1') \; , \; \vert \;\left.\right)}{4}
                                \hspace{-20pt} \treeop{\quad - \quad}{4} \hspace{-20pt}
                                \ingamma{\gamma_{3,2,1,1}\left(\right.\; (Q1') \; ,}{4}
                                \hspace{-5pt}\treeingamma{\\[-6pt]c}\hspace{-5pt}
                                \ingamma{, \; \vert \; , \; \vert \;\left.\right)}{4} \\[-22pt]
    \ingamma{(Q1')-(Q2'):}{4} & \ingamma{\gamma_{2,3,1}\left(\right.}{4}
                                \hspace{-5pt}\treeingamma{\\[-6pt]c}\hspace{-10pt}
                                \ingamma{, \; (Q2') \; , \; \vert \;\left.\right)}{4}
                                \hspace{-20pt} \treeop{\quad - \quad}{4} \hspace{-20pt}
                                \ingamma{\gamma_{3,2,1,1}\left(\right.\; (Q1') \; ,}{4}
                                \hspace{-5pt}\treeingamma{\\[-6pt]d}\hspace{-5pt}
                                \ingamma{, \; \vert \; , \; \vert \;\left.\right)}{4} \\[-22pt]
    \ingamma{(Q1')-(Q4'):}{4} & \ingamma{\gamma_{2,3,1}\left(\right.}{4}
                                \hspace{-5pt}\treeingamma{\\[-6pt]c}\hspace{-10pt}
                                \ingamma{, \; (Q4') \; , \; \vert \;\left.\right)}{4}
                                \hspace{-20pt} \treeop{\quad - \quad}{4} \hspace{-20pt}
                                \ingamma{\gamma_{3,2,1,1}\left(\right.\; (Q1') \; ,}{4}
                                \hspace{-5pt}\treeingamma{\\[-6pt]b}\hspace{-5pt}
                                \ingamma{, \; \vert \; , \; \vert \;\left.\right)}{4} \\[-22pt]
    \ingamma{(Q1')-(Q5'):}{4} & \ingamma{\gamma_{2,3,1}\left(\right.}{4}
                                \hspace{-5pt}\treeingamma{\\[-6pt]c}\hspace{-10pt}
                                \ingamma{, \; (Q5') \; , \; \vert \;\left.\right)}{4}
                                \hspace{-20pt} \treeop{\quad - \quad}{4} \hspace{-20pt}
                                \ingamma{\gamma_{3,2,1,1}\left(\right.\; (Q1') \; ,}{4}
                                \hspace{-5pt}\treeingamma{\\[-6pt]a}\hspace{-5pt}
                                \ingamma{, \; \vert \; , \; \vert \;\left.\right)}{4} \\[-22pt]
    \ingamma{(Q2')-(Q3'):}{4} & \ingamma{\gamma_{2,3,1}\left(\right.}{4}
                                \hspace{-5pt}\treeingamma{\\[-6pt]c}\hspace{-10pt}
                                \ingamma{, \; (Q3') \; , \; \vert \;\left.\right)}{4}
                                \hspace{-20pt} \treeop{\quad - \quad}{4} \hspace{-20pt}
                                \ingamma{\gamma_{3,2,1,1}\left(\right.\; (Q2') \; ,}{4}
                                \hspace{-5pt}\treeingamma{\\[-6pt]d}\hspace{-5pt}
                                \ingamma{, \; \vert \; , \; \vert \;\left.\right)}{4} \\[-22pt]
    \ingamma{(Q2')-(Q6'):}{4} & \ingamma{\gamma_{2,3,1}\left(\right.}{4}
                                \hspace{-5pt}\treeingamma{\\[-6pt]c}\hspace{-10pt}
                                \ingamma{, \; (Q6') \; , \; \vert \;\left.\right)}{4}
                                \hspace{-20pt} \treeop{\quad - \quad}{4} \hspace{-20pt}
                                \ingamma{\gamma_{3,2,1,1}\left(\right.\; (Q2') \; ,}{4}
                                \hspace{-5pt}\treeingamma{\\[-6pt]a}\hspace{-5pt}
                                \ingamma{, \; \vert \; , \; \vert \;\left.\right)}{4} \\[-22pt]
    \ingamma{(Q3')-(Q3'):}{4} & \ingamma{\gamma_{2,3,1}\left(\right.}{4}
                                \hspace{-5pt}\treeingamma{\\[-6pt]d}\hspace{-10pt}
                                \ingamma{, \; (Q3') \; , \; \vert \;\left.\right)}{4}
                                \hspace{-20pt} \treeop{\quad - \quad}{4} \hspace{-20pt}
                                \ingamma{\gamma_{3,2,1,1}\left(\right.\; (Q3') \; ,}{4}
                                \hspace{-5pt}\treeingamma{\\[-6pt]d}\hspace{-5pt}
                                \ingamma{, \; \vert \; , \; \vert \;\left.\right)}{4} \\[-22pt]
    \ingamma{(Q3')-(Q6'):}{4} & \ingamma{\gamma_{2,3,1}\left(\right.}{4}
                                \hspace{-5pt}\treeingamma{\\[-6pt]d}\hspace{-10pt}
                                \ingamma{, \; (Q6') \; , \; \vert \;\left.\right)}{4}
                                \hspace{-20pt} \treeop{\quad - \quad}{4} \hspace{-20pt}
                                \ingamma{\gamma_{3,2,1,1}\left(\right.\; (Q3') \; ,}{4}
                                \hspace{-5pt}\treeingamma{\\[-6pt]a}\hspace{-5pt}
                                \ingamma{, \; \vert \; , \; \vert \;\left.\right)}{4} \\[-22pt]
    \ingamma{(Q4')-(Q7'):}{4} & \ingamma{\gamma_{2,3,1}\left(\right.}{4}
                                \hspace{-5pt}\treeingamma{\\[-6pt]c}\hspace{-10pt}
                                \ingamma{, \; (Q7') \; , \; \vert \;\left.\right)}{4}
                                \hspace{-20pt} \treeop{\quad - \quad}{4} \hspace{-20pt}
                                \ingamma{\gamma_{3,2,1,1}\left(\right.\; (Q4') \; ,}{4}
                                \hspace{-5pt}\treeingamma{\\[-6pt]b}\hspace{-5pt}
                                \ingamma{, \; \vert \; , \; \vert \;\left.\right)}{4} \\[-22pt]
    \ingamma{(Q4')-(Q8'):}{4} & \ingamma{\gamma_{2,3,1}\left(\right.}{4}
                                \hspace{-5pt}\treeingamma{\\[-6pt]c}\hspace{-10pt}
                                \ingamma{, \; (Q8') \; , \; \vert \;\left.\right)}{4}
                                \hspace{-20pt} \treeop{\quad - \quad}{4} \hspace{-20pt}
                                \ingamma{\gamma_{3,2,1,1}\left(\right.\; (Q4') \; ,}{4}
                                \hspace{-5pt}\treeingamma{\\[-6pt]a}\hspace{-5pt}
                                \ingamma{, \; \vert \; , \; \vert \;\left.\right)}{4} \\[-22pt]
    \ingamma{(Q5')-(Q9'):}{4} & \ingamma{\gamma_{2,3,1}\left(\right.}{4}
                                \hspace{-5pt}\treeingamma{\\[-6pt]c}\hspace{-10pt}
                                \ingamma{, \; (Q9') \; , \; \vert \;\left.\right)}{4}
                                \hspace{-20pt} \treeop{\quad - \quad}{4} \hspace{-20pt}
                                \ingamma{\gamma_{3,2,1,1}\left(\right.\; (Q5') \; ,}{4}
                                \hspace{-5pt}\treeingamma{\\[-6pt]a}\hspace{-5pt}
                                \ingamma{, \; \vert \; , \; \vert \;\left.\right)}{4} \\[-22pt]
    \ingamma{(Q6')-(Q9'):}{4} & \ingamma{\gamma_{2,3,1}\left(\right.}{4}
                                \hspace{-5pt}\treeingamma{\\[-6pt]a}\hspace{-10pt}
                                \ingamma{, \; (Q9') \; , \; \vert \;\left.\right)}{4}
                                \hspace{-20pt} \treeop{\quad - \quad}{4} \hspace{-20pt}
                                \ingamma{\gamma_{3,2,1,1}\left(\right.\; (Q6') \; ,}{4}
                                \hspace{-5pt}\treeingamma{\\[-6pt]a}\hspace{-5pt}
                                \ingamma{, \; \vert \; , \; \vert \;\left.\right)}{4} \\[-22pt]
    \ingamma{(Q7')-(Q7'):}{4} & \ingamma{\gamma_{2,3,1}\left(\right.}{4}
                                \hspace{-5pt}\treeingamma{\\[-6pt]b}\hspace{-10pt}
                                \ingamma{, \; (Q7') \; , \; \vert \;\left.\right)}{4}
                                \hspace{-20pt} \treeop{\quad - \quad}{4} \hspace{-20pt}
                                \ingamma{\gamma_{3,2,1,1}\left(\right.\; (Q7') \; ,}{4}
                                \hspace{-5pt}\treeingamma{\\[-6pt]b}\hspace{-5pt}
                                \ingamma{, \; \vert \; , \; \vert \;\left.\right)}{4} \\[-22pt]
    \ingamma{(Q7')-(Q8'):}{4} & \ingamma{\gamma_{2,3,1}\left(\right.}{4}
                                \hspace{-5pt}\treeingamma{\\[-6pt]b}\hspace{-10pt}
                                \ingamma{, \; (Q8') \; , \; \vert \;\left.\right)}{4}
                                \hspace{-20pt} \treeop{\quad - \quad}{4} \hspace{-20pt}
                                \ingamma{\gamma_{3,2,1,1}\left(\right.\; (Q7') \; ,}{4}
                                \hspace{-5pt}\treeingamma{\\[-6pt]a}\hspace{-5pt}
                                \ingamma{, \; \vert \; , \; \vert \;\left.\right)}{4} \\[-22pt]
    \ingamma{(Q8')-(Q9'):}{4} & \ingamma{\gamma_{2,3,1}\left(\right.}{4}
                                \hspace{-5pt}\treeingamma{\\[-6pt]b}\hspace{-10pt}
                                \ingamma{, \; (Q9') \; , \; \vert \;\left.\right)}{4}
                                \hspace{-20pt} \treeop{\quad - \quad}{4} \hspace{-20pt}
                                \ingamma{\gamma_{3,2,1,1}\left(\right.\; (Q8') \; ,}{4}
                                \hspace{-5pt}\treeingamma{\\[-6pt]a}\hspace{-5pt}
                                \ingamma{, \; \vert \; , \; \vert \;\left.\right)}{4} \\[-22pt]
    \ingamma{(Q9')-(Q9'):}{4} & \ingamma{\gamma_{2,3,1}\left(\right.}{4}
                                \hspace{-5pt}\treeingamma{\\[-6pt]a}\hspace{-10pt}
                                \ingamma{, \; (Q9') \; , \; \vert \;\left.\right)}{4}
                                \hspace{-20pt} \treeop{\quad - \quad}{4} \hspace{-20pt}
                                \ingamma{\gamma_{3,2,1,1}\left(\right.\; (Q9') \; ,}{4}
                                \hspace{-5pt}\treeingamma{\\[-6pt]a}\hspace{-5pt}
                                \ingamma{, \; \vert \; , \; \vert \;\left.\right)}{4}
\end{align*} \\[-22pt]
We now reduce all these compositions using the elements in $S$. As before, we use subindexes to make clear which element
we use in each step of reduction. \\[-30pt]
\begin{align*}
    & \ingamma{\gamma_{2,3,1}\left(\right.}{4}
                                \hspace{-5pt}\treeingamma{\\[-6pt]c}\hspace{-10pt}
                                \ingamma{, \; (Q1') \; , \; \vert \;\left.\right)}{4}
                                \hspace{-20pt} \treeop{\quad - \quad}{4} \hspace{-20pt}
                                \ingamma{\gamma_{3,2,1,1}\left(\right.\; (Q1') \; ,}{4}
                                \hspace{-5pt}\treeingamma{\\[-6pt]c}\hspace{-5pt}
                                \ingamma{, \; \vert \; , \; \vert \;\left.\right)}{4}
    \\[-25pt]
    & \; \treeop{=}{25} \tiny
    \lltree{c}{c}{c}
    \treeop{-}{25}
    \lrtree{c}{c}{a}_{(Q1')}
    \treeop{-}{25}
    \lrtree{c}{c}{d}_{(Q1')}
    \treeop{-}{25}
    \lrtree{c}{c}{b}_{(Q1')}
    \treeop{-}{25}
    \lrtree{c}{c}{c}_{(Q1')}
    \\[-20pt]
    & \; \tiny
    \treeop{-}{25}
    \lltree{c}{c}{c}
    \treeop{+}{25}
    \corolla{c}{c}{a}_{(Q1')}
    \treeop{+}{25}
    \corolla{c}{c}{d}_{(Q1')}
    \treeop{+}{25}
    \corolla{c}{c}{b}_{(Q1')}
    \treeop{+}{25}
    \corolla{c}{c}{c}_{(Q1')}
    \\[-20pt]
    & \; \treeop{\equiv}{25} \tiny
    \treeop{-}{25}
    \rltree{c}{a}{a}_{(Q9')}
    \treeop{-}{25}
    \rltree{c}{d}{a}_{(Q6')}
    \treeop{-}{25}
    \rltree{c}{b}{a}_{(Q8')}
    \treeop{-}{25}
    \rltree{c}{c}{a}_{(Q5')}
    \treeop{-}{25}
    \rltree{c}{a}{d}
    \\[-20pt]
    & \; \tiny
    \treeop{-}{25}
    \rltree{c}{d}{d}_{(Q3')}
    \treeop{-}{25}
    \rltree{c}{b}{d}
    \treeop{-}{25}
    \rltree{c}{c}{d}_{(Q2')}
    \treeop{-}{25}
    \rltree{c}{a}{b}
    \treeop{-}{25}
    \rltree{c}{d}{b}
    \\[-20pt]
    & \; \tiny
    \treeop{-}{25}
    \rltree{c}{b}{b}_{(Q7')}
    \treeop{-}{25}
    \rltree{c}{c}{b}_{(Q4')}
    \treeop{-}{25}
    \rltree{c}{a}{c}
    \treeop{-}{25}
    \rltree{c}{d}{c}
    \treeop{-}{25}
    \rltree{c}{b}{c}
    \\[-20pt]
    & \; \tiny
    \treeop{-}{25}
    \rltree{c}{c}{c}_{(Q1')}
    \treeop{+}{25}
    \rrtree{c}{a}{a}
    \treeop{+}{25}
    \rrtree{c}{d}{a}
    \treeop{+}{25}
    \rrtree{c}{b}{a}
    \treeop{+}{25}
    \rrtree{c}{c}{a}
    \\[-20pt]
    & \; \tiny
    \treeop{+}{25}
    \rrtree{c}{a}{d}
    \treeop{+}{25}
    \rrtree{c}{d}{d}
    \treeop{+}{25}
    \rrtree{c}{b}{d}
    \treeop{+}{25}
    \rrtree{c}{c}{d}
    \treeop{+}{25}
    \rrtree{c}{a}{b}
    \\[-20pt]
    & \; \tiny
    \treeop{+}{25}
    \rrtree{c}{d}{b}
    \treeop{+}{25}
    \rrtree{c}{b}{b}
    \treeop{+}{25}
    \rrtree{c}{c}{b}
    \treeop{+}{25}
    \rrtree{c}{a}{c}
    \treeop{+}{25}
    \rrtree{c}{d}{c}
    \\[-20pt]
    & \; \tiny
    \treeop{+}{25}
    \rrtree{c}{b}{c}
    \treeop{+}{25}
    \rrtree{c}{c}{c}
    \\[-20pt]
    & \; \treeop{\equiv}{25} \tiny
    \rltree{c}{a}{d}
    \treeop{+}{25}
    \rltree{c}{a}{b}
    \treeop{+}{25}
    \rltree{c}{a}{c}
    \treeop{-}{25}
    \rrtree{c}{a}{a}
    \treeop{+}{25}
    \rltree{c}{d}{b}
    \treeop{-}{25}
    \rrtree{c}{a}{d}
    \\[-20pt]
    & \; \tiny
    \treeop{+}{25}
    \rltree{c}{b}{d}
    \treeop{-}{25}
    \rrtree{c}{a}{b}
    \treeop{-}{25}
    \rrtree{c}{a}{c}
    \treeop{-}{25}
    \rltree{c}{a}{d}
    \treeop{+}{25}
    \rltree{c}{d}{c}
    \treeop{-}{25}
    \rrtree{c}{d}{a}
    \\[-20pt]
    & \; \tiny
    \treeop{-}{25}
    \rrtree{c}{d}{d}
    \treeop{-}{25}
    \rltree{c}{b}{d}
    \treeop{-}{25}
    \rrtree{c}{d}{b}
    \treeop{-}{25}
    \rrtree{c}{d}{c}
    \treeop{-}{25}
    \rltree{c}{a}{b}
    \treeop{-}{25}
    \rltree{c}{d}{b}
    \\[-20pt]
    & \; \tiny
    \treeop{+}{25}
    \rltree{c}{b}{c}
    \treeop{-}{25}
    \rrtree{c}{b}{a}
    \treeop{-}{25}
    \rrtree{c}{b}{b}
    \treeop{-}{25}
    \rrtree{c}{b}{d}
    \treeop{-}{25}
    \rrtree{c}{b}{c}
    \treeop{-}{25}
    \rltree{c}{a}{c}
    \\[-20pt]
    & \; \tiny
    \treeop{-}{25}
    \rltree{c}{d}{c}
    \treeop{-}{25}
    \rltree{c}{b}{c}
    \treeop{-}{25}
    \rrtree{c}{c}{a}
    \treeop{-}{25}
    \rrtree{c}{c}{d}
    \treeop{-}{25}
    \rrtree{c}{c}{b}
    \treeop{-}{25}
    \rrtree{c}{c}{c}
    \\[-20pt]
    & \; \tiny
    \treeop{+}{25}
    \rrtree{c}{a}{a}
    \treeop{+}{25}
    \rrtree{c}{d}{a}
    \treeop{+}{25}
    \rrtree{c}{b}{a}
    \treeop{+}{25}
    \rrtree{c}{c}{a}
    \treeop{+}{25}
    \rrtree{c}{a}{d}
    \treeop{+}{25}
    \rrtree{c}{d}{d}
    \\[-20pt]
    & \; \tiny
    \treeop{+}{25}
    \rrtree{c}{b}{d}
    \treeop{+}{25}
    \rrtree{c}{c}{d}
    \treeop{+}{25}
    \rrtree{c}{a}{b}
    \treeop{+}{25}
    \rrtree{c}{d}{b}
    \treeop{+}{25}
    \rrtree{c}{b}{b}
    \treeop{+}{25}
    \rrtree{c}{c}{b}
    \\[-20pt]
    & \; \tiny
    \treeop{+}{25}
    \rrtree{c}{a}{c}
    \treeop{+}{25}
    \rrtree{c}{d}{c}
    \treeop{+}{25}
    \rrtree{c}{b}{c}
    \treeop{+}{25}
    \rrtree{c}{c}{c}
    \normalsize \treeop{= 0}{25}
    \\[-10pt]
    & \ingamma{\gamma_{2,3,1}\left(\right.}{4}
                                \hspace{-5pt}\treeingamma{\\[-6pt]c}\hspace{-10pt}
                                \ingamma{, \; (Q2') \; , \; \vert \;\left.\right)}{4}
                                \hspace{-20pt} \treeop{\quad - \quad}{4} \hspace{-20pt}
                                \ingamma{\gamma_{3,2,1,1}\left(\right.\; (Q1') \; ,}{4}
                                \hspace{-5pt}\treeingamma{\\[-6pt]d}\hspace{-5pt}
                                \ingamma{, \; \vert \; , \; \vert \;\left.\right)}{4}
    \\[-25pt]
    & \; \treeop{=}{25} \tiny
    \lltree{c}{c}{d}
    \treeop{-}{25}
    \lrtree{c}{d}{b}_{(Q2')}
    \treeop{-}{25}
    \lrtree{c}{d}{c}_{(Q2')}
    \treeop{-}{25}
    \lltree{c}{c}{d}
    \treeop{+}{25}
    \corolla{c}{d}{a}_{(Q2')}
    \\[-20pt]
    & \; \tiny
    \treeop{+}{25}
    \corolla{c}{d}{d}_{(Q2')}
    \treeop{+}{25}
    \corolla{c}{d}{b}_{(Q2')}
    \treeop{+}{25}
    \corolla{c}{d}{c}_{(Q2')}
    \\[-20pt]
    & \; \treeop{\equiv}{25} \tiny
    \treeop{-}{25}
    \rltree{d}{b}{b}_{(Q7')}
    \treeop{-}{25}
    \rltree{d}{c}{b}_{(Q4')}
    \treeop{-}{25}
    \rltree{d}{b}{c}
    \treeop{-}{25}
    \rltree{d}{c}{c}_{(Q1')}
    \treeop{+}{25}
    \rrtree{d}{b}{a}
    \\[-20pt]
    & \; \tiny
    \treeop{+}{25}
    \rrtree{d}{c}{a}
    \treeop{+}{25}
    \rrtree{d}{b}{d}
    \treeop{+}{25}
    \rrtree{d}{c}{d}
    \treeop{+}{25}
    \rrtree{d}{b}{b}
    \treeop{+}{25}
    \rrtree{d}{c}{b}
    \\[-20pt]
    & \; \tiny
    \treeop{+}{25}
    \rrtree{d}{b}{c}
    \treeop{+}{25}
    \rrtree{d}{c}{c}
    \\[-20pt]
    & \; \treeop{\equiv}{25} \tiny
    \rltree{d}{b}{c}
    \treeop{-}{25}
    \rrtree{d}{b}{a}
    \treeop{-}{25}
    \rrtree{d}{b}{b}
    \treeop{-}{25}
    \rrtree{d}{b}{d}
    \treeop{-}{25}
    \rrtree{d}{b}{c}
    \treeop{-}{25}
    \rltree{d}{b}{c}
    \\[-20pt]
    & \; \tiny
    \treeop{-}{25}
    \rrtree{d}{c}{a}
    \treeop{-}{25}
    \rrtree{d}{c}{d}
    \treeop{-}{25}
    \rrtree{d}{c}{b}
    \treeop{-}{25}
    \rrtree{d}{c}{c}
    \treeop{+}{25}
    \rrtree{d}{b}{a}
    \treeop{+}{25}
    \rrtree{d}{c}{a}
    \\[-20pt]
    & \; \tiny
    \treeop{+}{25}
    \rrtree{d}{b}{d}
    \treeop{+}{25}
    \rrtree{d}{c}{d}
    \treeop{+}{25}
    \rrtree{d}{b}{b}
    \treeop{+}{25}
    \rrtree{d}{c}{b}
    \treeop{+}{25}
    \rrtree{d}{b}{c}
    \treeop{+}{25}
    \rrtree{d}{c}{c}
    \normalsize \treeop{= 0}{25}
    \\[-10pt]
    & \ingamma{\gamma_{2,3,1}\left(\right.}{4}
                                \hspace{-5pt}\treeingamma{\\[-6pt]c}\hspace{-10pt}
                                \ingamma{, \; (Q4') \; , \; \vert \;\left.\right)}{4}
                                \hspace{-20pt} \treeop{\quad - \quad}{4} \hspace{-20pt}
                                \ingamma{\gamma_{3,2,1,1}\left(\right.\; (Q1') \; ,}{4}
                                \hspace{-5pt}\treeingamma{\\[-6pt]b}\hspace{-5pt}
                                \ingamma{, \; \vert \; , \; \vert \;\left.\right)}{4}
    \\[-25pt]
    & \; \treeop{=}{25} \tiny
    \lltree{c}{c}{b}
    \treeop{-}{25}
    \lrtree{c}{b}{d}_{(Q4')}
    \treeop{-}{25}
    \lrtree{c}{b}{c}_{(Q4')}
    \treeop{-}{25}
    \lltree{c}{c}{b}
    \treeop{+}{25}
    \corolla{c}{b}{a}_{(Q4')}
    \\[-20pt]
    & \; \tiny
    \treeop{+}{25}
    \corolla{c}{b}{d}_{(Q4')}
    \treeop{+}{25}
    \corolla{c}{b}{b}_{(Q4')}
    \treeop{+}{25}
    \corolla{c}{b}{c}_{(Q4')}
    \\[-20pt]
    & \; \treeop{\equiv}{25} \tiny
    \treeop{-}{25}
    \rltree{b}{d}{d}_{(Q3')}
    \treeop{-}{25}
    \rltree{b}{c}{d}
    \treeop{-}{25}
    \rltree{b}{d}{c}_{(Q2')}
    \treeop{-}{25}
    \rltree{b}{c}{c}_{(Q1')}
    \treeop{+}{25}
    \rrtree{b}{d}{a}
    \\[-20pt]
    & \; \tiny
    \treeop{+}{25}
    \rrtree{b}{c}{a}
    \treeop{+}{25}
    \rrtree{b}{d}{d}
    \treeop{+}{25}
    \rrtree{b}{c}{d}
    \treeop{+}{25}
    \rrtree{b}{d}{b}
    \treeop{+}{25}
    \rrtree{b}{c}{b}
    \\[-20pt]
    & \; \tiny
    \treeop{+}{25}
    \rrtree{b}{d}{c}
    \treeop{+}{25}
    \rrtree{b}{c}{c}
    \\[-20pt]
    & \; \treeop{\equiv}{25} \tiny
    \rltree{b}{d}{c}
    \treeop{-}{25}
    \rrtree{b}{d}{a}
    \treeop{-}{25}
    \rrtree{b}{d}{d}
    \treeop{-}{25}
    \rrtree{b}{d}{b}
    \treeop{-}{25}
    \rrtree{b}{d}{c}
    \treeop{-}{25}
    \rltree{b}{d}{c}
    \\[-20pt]
    & \; \tiny
    \treeop{-}{25}
    \rrtree{b}{c}{a}
    \treeop{-}{25}
    \rrtree{b}{c}{d}
    \treeop{-}{25}
    \rrtree{b}{c}{b}
    \treeop{-}{25}
    \rrtree{b}{c}{c}
    \treeop{+}{25}
    \rrtree{b}{d}{a}
    \treeop{+}{25}
    \rrtree{b}{c}{a}
    \\[-20pt]
    & \; \tiny
    \treeop{+}{25}
    \rrtree{b}{d}{d}
    \treeop{+}{25}
    \rrtree{b}{c}{d}
    \treeop{+}{25}
    \rrtree{b}{d}{b}
    \treeop{+}{25}
    \rrtree{b}{c}{b}
    \treeop{+}{25}
    \rrtree{b}{d}{c}
    \treeop{+}{25}
    \rrtree{b}{c}{c}
    \normalsize \treeop{= 0}{25}
    \\[-10pt]
    & \ingamma{\gamma_{2,3,1}\left(\right.}{4}
                                \hspace{-5pt}\treeingamma{\\[-6pt]c}\hspace{-10pt}
                                \ingamma{, \; (Q5') \; , \; \vert \;\left.\right)}{4}
                                \hspace{-20pt} \treeop{\quad - \quad}{4} \hspace{-20pt}
                                \ingamma{\gamma_{3,2,1,1}\left(\right.\; (Q1') \; ,}{4}
                                \hspace{-5pt}\treeingamma{\\[-6pt]a}\hspace{-5pt}
                                \ingamma{, \; \vert \; , \; \vert \;\left.\right)}{4}
    \\[-25pt]
    & \; \treeop{=}{25} \tiny
    \lltree{c}{c}{a}
    \treeop{-}{25}
    \lrtree{c}{a}{c}_{(Q5')}
    \treeop{-}{25}
    \lltree{c}{c}{a}
    \treeop{+}{25}
    \corolla{c}{a}{a}_{(Q5')}
    \treeop{+}{25}
    \corolla{c}{a}{d}_{(Q5')}
    \\[-20pt]
    & \; \tiny
    \treeop{+}{25}
    \corolla{c}{a}{b}_{(Q5')}
    \treeop{+}{25}
    \corolla{c}{a}{c}_{(Q5')}
    \\[-20pt]
    & \; \treeop{\equiv}{25} \tiny
    \treeop{-}{25}
    \rltree{a}{c}{c}_{(Q1')}
    \treeop{+}{25}
    \rrtree{a}{c}{a}
    \treeop{+}{25}
    \rrtree{a}{c}{d}
    \treeop{+}{25}
    \rrtree{a}{c}{b}
    \treeop{+}{25}
    \rrtree{a}{c}{c}
    \\[-20pt]
    & \; \tiny
    \\[-20pt]
    & \; \treeop{\equiv}{25} \tiny
    \treeop{-}{25}
    \rrtree{a}{c}{a}
    \treeop{-}{25}
    \rrtree{a}{c}{d}
    \treeop{-}{25}
    \rrtree{a}{c}{b}
    \treeop{-}{25}
    \rrtree{a}{c}{c}
    \treeop{+}{25}
    \rrtree{a}{c}{a}
    \treeop{+}{25}
    \rrtree{a}{c}{d}
    \\[-20pt]
    & \; \tiny
    \treeop{+}{25}
    \rrtree{a}{c}{b}
    \treeop{+}{25}
    \rrtree{a}{c}{c}
    \normalsize \treeop{= 0}{25}
    \\[-10pt]
    & \ingamma{\gamma_{2,3,1}\left(\right.}{4}
                                \hspace{-5pt}\treeingamma{\\[-6pt]c}\hspace{-10pt}
                                \ingamma{, \; (Q3') \; , \; \vert \;\left.\right)}{4}
                                \hspace{-20pt} \treeop{\quad - \quad}{4} \hspace{-20pt}
                                \ingamma{\gamma_{3,2,1,1}\left(\right.\; (Q2') \; ,}{4}
                                \hspace{-5pt}\treeingamma{\\[-6pt]d}\hspace{-5pt}
                                \ingamma{, \; \vert \; , \; \vert \;\left.\right)}{4}
    \\[-25pt]
    & \; \treeop{=}{25} \tiny
    \lltree{c}{d}{d}
    \treeop{+}{25}
    \lltree{c}{d}{c}_{(Q2')}
    \treeop{-}{25}
    \lrtree{c}{d}{a}_{(Q2')}
    \treeop{-}{25}
    \lrtree{c}{d}{d}_{(Q2')}
    \treeop{-}{25}
    \lltree{c}{d}{d}
    \\[-20pt]
    & \; \tiny
    \treeop{+}{25}
    \corolla{d}{d}{b}_{(Q3')}
    \treeop{+}{25}
    \corolla{d}{d}{c}_{(Q3')}
    \\[-20pt]
    & \; \treeop{\equiv}{25} \tiny
    \corolla{d}{c}{b}
    \treeop{+}{25}
    \corolla{d}{c}{c}
    \treeop{-}{25}
    \rltree{d}{b}{a}_{(Q8')}
    \treeop{-}{25}
    \rltree{d}{c}{a}_{(Q5')}
    \treeop{-}{25}
    \rltree{d}{b}{d}
    \\[-20pt]
    & \; \tiny
    \treeop{-}{25}
    \rltree{d}{c}{d}_{(Q2')}
    \treeop{-}{25}
    \corolla{d}{c}{b}
    \treeop{+}{25}
    \rrtree{d}{a}{b}
    \treeop{+}{25}
    \rrtree{d}{d}{b}
    \treeop{-}{25}
    \corolla{d}{c}{c}
    \\[-20pt]
    & \; \tiny
    \treeop{+}{25}
    \rrtree{d}{a}{c}
    \treeop{+}{25}
    \rrtree{d}{d}{c}
    \\[-20pt]
    & \; \treeop{\equiv}{25} \tiny
    \rltree{d}{b}{d}
    \treeop{-}{25}
    \rrtree{d}{a}{b}
    \treeop{-}{25}
    \rrtree{d}{a}{c}
    \treeop{-}{25}
    \rltree{d}{b}{d}
    \treeop{-}{25}
    \rrtree{d}{d}{b}
    \treeop{-}{25}
    \rrtree{d}{d}{c}
    \\[-20pt]
    & \; \tiny
    \treeop{+}{25}
    \rrtree{d}{a}{b}
    \treeop{+}{25}
    \rrtree{d}{d}{b}
    \treeop{+}{25}
    \rrtree{d}{a}{c}
    \treeop{+}{25}
    \rrtree{d}{d}{c}
    \normalsize \treeop{= 0}{25}
    \\[-10pt]
    & \ingamma{\gamma_{2,3,1}\left(\right.}{4}
                                \hspace{-5pt}\treeingamma{\\[-6pt]c}\hspace{-10pt}
                                \ingamma{, \; (Q6') \; , \; \vert \;\left.\right)}{4}
                                \hspace{-20pt} \treeop{\quad - \quad}{4} \hspace{-20pt}
                                \ingamma{\gamma_{3,2,1,1}\left(\right.\; (Q2') \; ,}{4}
                                \hspace{-5pt}\treeingamma{\\[-6pt]a}\hspace{-5pt}
                                \ingamma{, \; \vert \; , \; \vert \;\left.\right)}{4}
    \\[-25pt]
    & \; \treeop{=}{25} \tiny
    \lltree{c}{d}{a}
    \treeop{+}{25}
    \lltree{c}{d}{b}_{(Q2')}
    \treeop{-}{25}
    \lrtree{c}{a}{d}_{(Q5')}
    \treeop{-}{25}
    \lltree{c}{d}{a}
    \treeop{+}{25}
    \corolla{d}{a}{b}_{(Q6')}
    \\[-20pt]
    & \; \tiny
    \treeop{+}{25}
    \corolla{d}{a}{c}_{(Q6')}
    \\[-20pt]
    & \; \treeop{\equiv}{25} \tiny
    \corolla{d}{b}{b}
    \treeop{+}{25}
    \corolla{d}{b}{c}
    \treeop{-}{25}
    \rltree{a}{c}{d}_{(Q2')}
    \treeop{-}{25}
    \corolla{d}{b}{b}
    \treeop{+}{25}
    \rrtree{a}{d}{b}
    \\[-20pt]
    & \; \tiny
    \treeop{-}{25}
    \corolla{d}{b}{c}
    \treeop{+}{25}
    \rrtree{a}{d}{c}
    \\[-20pt]
    & \; \treeop{\equiv}{25} \tiny
    \treeop{-}{25}
    \rrtree{a}{d}{b}
    \treeop{-}{25}
    \rrtree{a}{d}{c}
    \treeop{+}{25}
    \rrtree{a}{d}{b}
    \treeop{+}{25}
    \rrtree{a}{d}{c}
    \normalsize \treeop{= 0}{25}
    \\[-10pt]
    & \ingamma{\gamma_{2,3,1}\left(\right.}{4}
                                \hspace{-5pt}\treeingamma{\\[-6pt]d}\hspace{-10pt}
                                \ingamma{, \; (Q3') \; , \; \vert \;\left.\right)}{4}
                                \hspace{-20pt} \treeop{\quad - \quad}{4} \hspace{-20pt}
                                \ingamma{\gamma_{3,2,1,1}\left(\right.\; (Q3') \; ,}{4}
                                \hspace{-5pt}\treeingamma{\\[-6pt]d}\hspace{-5pt}
                                \ingamma{, \; \vert \; , \; \vert \;\left.\right)}{4}
    \\[-25pt]
    & \; \treeop{=}{25} \tiny
    \lltree{d}{d}{d}
    \treeop{+}{25}
    \lltree{d}{d}{c}_{(Q3')}
    \treeop{-}{25}
    \lrtree{d}{d}{a}_{(Q3')}
    \treeop{-}{25}
    \lrtree{d}{d}{d}_{(Q3')}
    \treeop{-}{25}
    \lltree{d}{d}{d}
    \\[-20pt]
    & \; \tiny
    \treeop{-}{25}
    \lltree{d}{c}{d}_{(Q3')}
    \treeop{+}{25}
    \corolla{d}{d}{a}_{(Q3')}
    \treeop{+}{25}
    \corolla{d}{d}{d}_{(Q3')}
    \\[-20pt]
    & \; \treeop{\equiv}{25} \tiny
    \treeop{-}{25}
    \lltree{d}{c}{c}_{(Q1')}
    \treeop{+}{25}
    \corolla{d}{c}{a}
    \treeop{+}{25}
    \corolla{d}{c}{d}
    \treeop{+}{25}
    \lrtree{d}{c}{a}
    \treeop{-}{25}
    \rltree{d}{a}{a}_{(Q9')}
    \\[-20pt]
    & \; \tiny
    \treeop{-}{25}
    \rltree{d}{d}{a}_{(Q6')}
    \treeop{+}{25}
    \lrtree{d}{c}{d}
    \treeop{-}{25}
    \rltree{d}{a}{d}
    \treeop{-}{25}
    \rltree{d}{d}{d}_{(Q3')}
    \treeop{-}{25}
    \lrtree{d}{d}{b}_{(Q3')}
    \\[-20pt]
    & \; \tiny
    \treeop{-}{25}
    \lrtree{d}{d}{c}_{(Q3')}
    \treeop{-}{25}
    \corolla{d}{c}{a}
    \treeop{+}{25}
    \rrtree{d}{a}{a}
    \treeop{+}{25}
    \rrtree{d}{d}{a}
    \treeop{-}{25}
    \corolla{d}{c}{d}
    \\[-20pt]
    & \; \tiny
    \treeop{+}{25}
    \rrtree{d}{a}{d}
    \treeop{+}{25}
    \rrtree{d}{d}{d}
    \\[-20pt]
    & \; \treeop{\equiv}{25} \tiny
    \treeop{-}{25}
    \lrtree{d}{c}{a}
    \treeop{-}{25}
    \lrtree{d}{c}{d}
    \treeop{-}{25}
    \lrtree{d}{c}{b}
    \treeop{-}{25}
    \lrtree{d}{c}{c}
    \treeop{+}{25}
    \lrtree{d}{c}{a}
    \treeop{+}{25}
    \rltree{d}{a}{d}
    \\[-20pt]
    & \; \tiny
    \treeop{+}{25}
    \rltree{d}{a}{b}
    \treeop{+}{25}
    \rltree{d}{a}{c}
    \treeop{-}{25}
    \rrtree{d}{a}{a}
    \treeop{+}{25}
    \rltree{d}{d}{b}
    \treeop{-}{25}
    \rrtree{d}{a}{d}
    \treeop{+}{25}
    \lrtree{d}{c}{d}
    \\[-20pt]
    & \; \tiny
    \treeop{-}{25}
    \rltree{d}{a}{d}
    \treeop{+}{25}
    \rltree{d}{d}{c}
    \treeop{-}{25}
    \rrtree{d}{d}{a}
    \treeop{-}{25}
    \rrtree{d}{d}{d}
    \treeop{+}{25}
    \lrtree{d}{c}{b}
    \treeop{-}{25}
    \rltree{d}{a}{b}
    \\[-20pt]
    & \; \tiny
    \treeop{-}{25}
    \rltree{d}{d}{b}
    \treeop{+}{25}
    \lrtree{d}{c}{c}
    \treeop{-}{25}
    \rltree{d}{a}{c}
    \treeop{-}{25}
    \rltree{d}{d}{c}
    \treeop{+}{25}
    \rrtree{d}{a}{a}
    \treeop{+}{25}
    \rrtree{d}{d}{a}
    \\[-20pt]
    & \; \tiny
    \treeop{+}{25}
    \rrtree{d}{a}{d}
    \treeop{+}{25}
    \rrtree{d}{d}{d}
    \normalsize \treeop{= 0}{25}
    \\[-10pt]
    & \ingamma{\gamma_{2,3,1}\left(\right.}{4}
                                \hspace{-5pt}\treeingamma{\\[-6pt]d}\hspace{-10pt}
                                \ingamma{, \; (Q6') \; , \; \vert \;\left.\right)}{4}
                                \hspace{-20pt} \treeop{\quad - \quad}{4} \hspace{-20pt}
                                \ingamma{\gamma_{3,2,1,1}\left(\right.\; (Q3') \; ,}{4}
                                \hspace{-5pt}\treeingamma{\\[-6pt]a}\hspace{-5pt}
                                \ingamma{, \; \vert \; , \; \vert \;\left.\right)}{4}
    \\[-25pt]
    & \; \treeop{=}{25} \tiny
    \lltree{d}{d}{a}
    \treeop{+}{25}
    \lltree{d}{d}{b}_{(Q3')}
    \treeop{-}{25}
    \lrtree{d}{a}{d}_{(Q6')}
    \treeop{-}{25}
    \lltree{d}{d}{a}
    \treeop{-}{25}
    \lltree{d}{c}{a}
    \\[-20pt]
    & \; \tiny
    \treeop{+}{25}
    \corolla{d}{a}{a}_{(Q6')}
    \treeop{+}{25}
    \corolla{d}{a}{d}_{(Q6')}
    \\[-20pt]
    & \; \treeop{\equiv}{25} \tiny
    \treeop{-}{25}
    \lltree{d}{c}{b}_{(Q4')}
    \treeop{+}{25}
    \corolla{d}{b}{a}
    \treeop{+}{25}
    \corolla{d}{b}{d}
    \treeop{+}{25}
    \lrtree{d}{b}{d}
    \treeop{-}{25}
    \rltree{a}{d}{d}_{(Q3')}
    \\[-20pt]
    & \; \tiny
    \treeop{-}{25}
    \lrtree{d}{a}{c}_{(Q6')}
    \treeop{-}{25}
    \corolla{d}{b}{a}
    \treeop{+}{25}
    \rrtree{a}{d}{a}
    \treeop{-}{25}
    \corolla{d}{b}{d}
    \treeop{+}{25}
    \rrtree{a}{d}{d}
    \\[-20pt]
    & \; \treeop{\equiv}{25} \tiny
    \treeop{-}{25}
    \lrtree{d}{b}{d}
    \treeop{-}{25}
    \lrtree{d}{b}{c}
    \treeop{+}{25}
    \lrtree{d}{b}{d}
    \treeop{+}{25}
    \rltree{a}{d}{c}
    \treeop{-}{25}
    \rrtree{a}{d}{a}
    \treeop{-}{25}
    \rrtree{a}{d}{d}
    \\[-20pt]
    & \; \tiny
    \treeop{+}{25}
    \lrtree{d}{b}{c}
    \treeop{-}{25}
    \rltree{a}{d}{c}
    \treeop{+}{25}
    \rrtree{a}{d}{a}
    \treeop{+}{25}
    \rrtree{a}{d}{d}
    \normalsize \treeop{= 0}{25}
    \\[-10pt]
    & \ingamma{\gamma_{2,3,1}\left(\right.}{4}
                                \hspace{-5pt}\treeingamma{\\[-6pt]c}\hspace{-10pt}
                                \ingamma{, \; (Q7') \; , \; \vert \;\left.\right)}{4}
                                \hspace{-20pt} \treeop{\quad - \quad}{4} \hspace{-20pt}
                                \ingamma{\gamma_{3,2,1,1}\left(\right.\; (Q4') \; ,}{4}
                                \hspace{-5pt}\treeingamma{\\[-6pt]b}\hspace{-5pt}
                                \ingamma{, \; \vert \; , \; \vert \;\left.\right)}{4}
    \\[-25pt]
    & \; \treeop{=}{25} \tiny
    \lltree{c}{b}{b}
    \treeop{+}{25}
    \lltree{c}{b}{c}_{(Q4')}
    \treeop{-}{25}
    \lrtree{c}{b}{a}_{(Q4')}
    \treeop{-}{25}
    \lrtree{c}{b}{b}_{(Q4')}
    \treeop{-}{25}
    \lltree{c}{b}{b}
    \\[-20pt]
    & \; \tiny
    \treeop{+}{25}
    \corolla{b}{b}{d}_{(Q7')}
    \treeop{+}{25}
    \corolla{b}{b}{c}_{(Q7')}
    \\[-20pt]
    & \; \treeop{\equiv}{25} \tiny
    \corolla{b}{c}{d}
    \treeop{+}{25}
    \corolla{b}{c}{c}
    \treeop{-}{25}
    \rltree{b}{d}{a}_{(Q6')}
    \treeop{-}{25}
    \rltree{b}{c}{a}_{(Q5')}
    \treeop{-}{25}
    \rltree{b}{d}{b}
    \\[-20pt]
    & \; \tiny
    \treeop{-}{25}
    \rltree{b}{c}{b}_{(Q4')}
    \treeop{-}{25}
    \corolla{b}{c}{d}
    \treeop{+}{25}
    \rrtree{b}{a}{d}
    \treeop{+}{25}
    \rrtree{b}{b}{d}
    \treeop{-}{25}
    \corolla{b}{c}{c}
    \\[-20pt]
    & \; \tiny
    \treeop{+}{25}
    \rrtree{b}{a}{c}
    \treeop{+}{25}
    \rrtree{b}{b}{c}
    \\[-20pt]
    & \; \treeop{\equiv}{25} \tiny
    \rltree{b}{d}{b}
    \treeop{-}{25}
    \rrtree{b}{a}{d}
    \treeop{-}{25}
    \rrtree{b}{a}{c}
    \treeop{-}{25}
    \rltree{b}{d}{b}
    \treeop{-}{25}
    \rrtree{b}{b}{d}
    \treeop{-}{25}
    \rrtree{b}{b}{c}
    \\[-20pt]
    & \; \tiny
    \treeop{+}{25}
    \rrtree{b}{a}{d}
    \treeop{+}{25}
    \rrtree{b}{b}{d}
    \treeop{+}{25}
    \rrtree{b}{a}{c}
    \treeop{+}{25}
    \rrtree{b}{b}{c}
    \normalsize \treeop{= 0}{25}
    \\[-10pt]
    & \ingamma{\gamma_{2,3,1}\left(\right.}{4}
                                \hspace{-5pt}\treeingamma{\\[-6pt]c}\hspace{-10pt}
                                \ingamma{, \; (Q8') \; , \; \vert \;\left.\right)}{4}
                                \hspace{-20pt} \treeop{\quad - \quad}{4} \hspace{-20pt}
                                \ingamma{\gamma_{3,2,1,1}\left(\right.\; (Q4') \; ,}{4}
                                \hspace{-5pt}\treeingamma{\\[-6pt]a}\hspace{-5pt}
                                \ingamma{, \; \vert \; , \; \vert \;\left.\right)}{4}
    \\[-25pt]
    & \; \treeop{=}{25} \tiny
    \lltree{c}{b}{a}
    \treeop{+}{25}
    \lltree{c}{b}{d}_{(Q4')}
    \treeop{-}{25}
    \lrtree{c}{a}{b}_{(Q5')}
    \treeop{-}{25}
    \lltree{c}{b}{a}
    \treeop{+}{25}
    \corolla{b}{a}{d}_{(Q8')}
    \\[-20pt]
    & \; \tiny
    \treeop{+}{25}
    \corolla{b}{a}{c}_{(Q8')}
    \\[-20pt]
    & \; \treeop{\equiv}{25} \tiny
    \corolla{b}{d}{d}
    \treeop{+}{25}
    \corolla{b}{d}{c}
    \treeop{-}{25}
    \rltree{a}{c}{b}_{(Q4')}
    \treeop{-}{25}
    \corolla{b}{d}{d}
    \treeop{+}{25}
    \rrtree{a}{b}{d}
    \\[-20pt]
    & \; \tiny
    \treeop{-}{25}
    \corolla{b}{d}{c}
    \treeop{+}{25}
    \rrtree{a}{b}{c}
    \\[-20pt]
    & \; \treeop{\equiv}{25} \tiny
    \treeop{-}{25}
    \rrtree{a}{b}{d}
    \treeop{-}{25}
    \rrtree{a}{b}{c}
    \treeop{+}{25}
    \rrtree{a}{b}{d}
    \treeop{+}{25}
    \rrtree{a}{b}{c}
    \normalsize \treeop{= 0}{25}
    \\[-10pt]
    & \ingamma{\gamma_{2,3,1}\left(\right.}{4}
                                \hspace{-5pt}\treeingamma{\\[-6pt]c}\hspace{-10pt}
                                \ingamma{, \; (Q9') \; , \; \vert \;\left.\right)}{4}
                                \hspace{-20pt} \treeop{\quad - \quad}{4} \hspace{-20pt}
                                \ingamma{\gamma_{3,2,1,1}\left(\right.\; (Q5') \; ,}{4}
                                \hspace{-5pt}\treeingamma{\\[-6pt]a}\hspace{-5pt}
                                \ingamma{, \; \vert \; , \; \vert \;\left.\right)}{4}
    \\[-25pt]
    & \; \treeop{=}{25} \tiny
    \lltree{c}{a}{a}
    \treeop{+}{25}
    \lltree{c}{a}{d}_{(Q5')}
    \treeop{+}{25}
    \lltree{c}{a}{b}_{(Q5')}
    \treeop{+}{25}
    \lltree{c}{a}{c}_{(Q5')}
    \treeop{-}{25}
    \lrtree{c}{a}{a}_{(Q5')}
    \\[-20pt]
    & \; \tiny
    \treeop{-}{25}
    \lltree{c}{a}{a}
    \treeop{+}{25}
    \corolla{a}{a}{c}_{(Q9')}
    \\[-20pt]
    & \; \treeop{\equiv}{25} \tiny
    \corolla{a}{d}{c}
    \treeop{+}{25}
    \corolla{a}{b}{c}
    \treeop{+}{25}
    \corolla{a}{c}{c}
    \treeop{-}{25}
    \rltree{a}{c}{a}_{(Q5')}
    \treeop{-}{25}
    \corolla{a}{d}{c}
    \treeop{-}{25}
    \corolla{a}{b}{c}
    \\[-20pt]
    & \; \tiny
    \treeop{-}{25}
    \corolla{a}{c}{c}
    \treeop{+}{25}
    \rrtree{a}{a}{c}
    \; \treeop{\equiv}{25} \tiny
    \treeop{-}{25}
    \rrtree{a}{a}{c}
    \treeop{+}{25}
    \rrtree{a}{a}{c}
    \normalsize \treeop{= 0}{25}
    \\[-10pt]
    & \ingamma{\gamma_{2,3,1}\left(\right.}{4}
                                \hspace{-5pt}\treeingamma{\\[-6pt]a}\hspace{-10pt}
                                \ingamma{, \; (Q9') \; , \; \vert \;\left.\right)}{4}
                                \hspace{-20pt} \treeop{\quad - \quad}{4} \hspace{-20pt}
                                \ingamma{\gamma_{3,2,1,1}\left(\right.\; (Q6') \; ,}{4}
                                \hspace{-5pt}\treeingamma{\\[-6pt]a}\hspace{-5pt}
                                \ingamma{, \; \vert \; , \; \vert \;\left.\right)}{4}
    \\[-25pt]
    & \; \treeop{=}{25} \tiny
    \lltree{d}{a}{a}
    \treeop{+}{25}
    \lltree{d}{a}{d}_{(Q6')}
    \treeop{+}{25}
    \lltree{d}{a}{b}_{(Q6')}
    \treeop{+}{25}
    \lltree{d}{a}{c}_{(Q6')}
    \treeop{-}{25}
    \lrtree{d}{a}{a}_{(Q6')}
    \\[-20pt]
    & \; \tiny
    \treeop{-}{25}
    \lltree{d}{a}{a}
    \treeop{-}{25}
    \lltree{d}{b}{a}_{(Q8')}
    \treeop{+}{25}
    \corolla{a}{a}{d}_{(Q9')}
    \\[-20pt]
    & \; \treeop{\equiv}{25} \tiny
    \treeop{-}{25}
    \lltree{d}{b}{d}
    \treeop{+}{25}
    \corolla{a}{d}{d}
    \treeop{-}{25}
    \lltree{d}{b}{b}_{(Q7')}
    \treeop{+}{25}
    \corolla{a}{b}{d}
    \treeop{-}{25}
    \lltree{d}{b}{c}
    \\[-20pt]
    & \; \tiny
    \treeop{+}{25}
    \corolla{a}{c}{d}
    \treeop{+}{25}
    \lrtree{d}{b}{a}
    \treeop{-}{25}
    \rltree{a}{d}{a}_{(Q6')}
    \treeop{+}{25}
    \lltree{d}{b}{d}
    \treeop{-}{25}
    \lrtree{d}{a}{b}_{(Q6')}
    \\[-20pt]
    & \; \tiny
    \treeop{-}{25}
    \corolla{a}{d}{d}
    \treeop{-}{25}
    \corolla{a}{b}{d}
    \treeop{-}{25}
    \corolla{a}{c}{d}
    \treeop{+}{25}
    \rrtree{a}{a}{d}
    \\[-20pt]
    & \; \treeop{\equiv}{25} \tiny
    \lltree{d}{b}{c}
    \treeop{-}{25}
    \lrtree{d}{b}{a}
    \treeop{-}{25}
    \lrtree{d}{b}{b}
    \treeop{-}{25}
    \lltree{d}{b}{c}
    \treeop{+}{25}
    \lrtree{d}{b}{a}
    \treeop{+}{25}
    \rltree{a}{d}{b}
    \\[-20pt]
    & \; \tiny
    \treeop{-}{25}
    \rrtree{a}{a}{d}
    \treeop{+}{25}
    \lrtree{d}{b}{b}
    \treeop{-}{25}
    \rltree{a}{d}{b}
    \treeop{+}{25}
    \rrtree{a}{a}{d}
    \normalsize \treeop{= 0}{25}
    \\[-10pt]
    & \ingamma{\gamma_{2,3,1}\left(\right.}{4}
                                \hspace{-5pt}\treeingamma{\\[-6pt]b}\hspace{-10pt}
                                \ingamma{, \; (Q7') \; , \; \vert \;\left.\right)}{4}
                                \hspace{-20pt} \treeop{\quad - \quad}{4} \hspace{-20pt}
                                \ingamma{\gamma_{3,2,1,1}\left(\right.\; (Q7') \; ,}{4}
                                \hspace{-5pt}\treeingamma{\\[-6pt]b}\hspace{-5pt}
                                \ingamma{, \; \vert \; , \; \vert \;\left.\right)}{4}
    \\[-25pt]
    & \; \treeop{=}{25} \tiny
    \lltree{b}{b}{b}
    \treeop{+}{25}
    \lltree{b}{b}{c}_{(Q7')}
    \treeop{-}{25}
    \lrtree{b}{b}{a}_{(Q7')}
    \treeop{-}{25}
    \lrtree{b}{b}{b}_{(Q7')}
    \treeop{-}{25}
    \lltree{b}{b}{b}
    \\[-20pt]
    & \; \tiny
    \treeop{-}{25}
    \lltree{b}{c}{b}_{(Q4')}
    \treeop{+}{25}
    \corolla{b}{b}{a}_{(Q7')}
    \treeop{+}{25}
    \corolla{b}{b}{b}_{(Q7')}
    \\[-20pt]
    & \; \treeop{\equiv}{25} \tiny
    \treeop{-}{25}
    \lltree{b}{c}{c}_{(Q1')}
    \treeop{+}{25}
    \corolla{b}{c}{a}
    \treeop{+}{25}
    \corolla{b}{c}{b}
    \treeop{+}{25}
    \lrtree{b}{c}{a}
    \treeop{-}{25}
    \rltree{b}{a}{a}_{(Q9')}
    \\[-20pt]
    & \; \tiny
    \treeop{-}{25}
    \rltree{b}{b}{a}_{(Q8')}
    \treeop{+}{25}
    \lrtree{b}{c}{b}
    \treeop{-}{25}
    \rltree{b}{a}{b}
    \treeop{-}{25}
    \rltree{b}{b}{b}_{(Q7')}
    \treeop{-}{25}
    \lrtree{b}{b}{d}_{(Q7')}
    \\[-20pt]
    & \; \tiny
    \treeop{-}{25}
    \lrtree{b}{b}{c}_{(Q7')}
    \treeop{-}{25}
    \corolla{b}{c}{a}
    \treeop{+}{25}
    \rrtree{b}{a}{a}
    \treeop{+}{25}
    \rrtree{b}{b}{a}
    \treeop{-}{25}
    \corolla{b}{c}{b}
    \\[-20pt]
    & \; \tiny
    \treeop{+}{25}
    \rrtree{b}{a}{b}
    \treeop{+}{25}
    \rrtree{b}{b}{b}
    \\[-20pt]
    & \; \treeop{\equiv}{25} \tiny
    \treeop{-}{25}
    \lrtree{b}{c}{a}
    \treeop{-}{25}
    \lrtree{b}{c}{d}
    \treeop{-}{25}
    \lrtree{b}{c}{b}
    \treeop{-}{25}
    \lrtree{b}{c}{c}
    \treeop{+}{25}
    \lrtree{b}{c}{a}
    \treeop{+}{25}
    \rltree{b}{a}{d}
    \\[-20pt]
    & \; \tiny
    \treeop{+}{25}
    \rltree{b}{a}{b}
    \treeop{+}{25}
    \rltree{b}{a}{c}
    \treeop{-}{25}
    \rrtree{b}{a}{a}
    \treeop{+}{25}
    \rltree{b}{b}{d}
    \treeop{-}{25}
    \rrtree{b}{a}{b}
    \treeop{+}{25}
    \lrtree{b}{c}{b}
    \\[-20pt]
    & \; \tiny
    \treeop{-}{25}
    \rltree{b}{a}{b}
    \treeop{+}{25}
    \rltree{b}{b}{c}
    \treeop{-}{25}
    \rrtree{b}{b}{a}
    \treeop{-}{25}
    \rrtree{b}{b}{b}
    \treeop{+}{25}
    \lrtree{b}{c}{d}
    \treeop{-}{25}
    \rltree{b}{a}{d}
    \\[-20pt]
    & \; \tiny
    \treeop{-}{25}
    \rltree{b}{b}{d}
    \treeop{+}{25}
    \lrtree{b}{c}{c}
    \treeop{-}{25}
    \rltree{b}{a}{c}
    \treeop{-}{25}
    \rltree{b}{b}{c}
    \treeop{+}{25}
    \rrtree{b}{a}{a}
    \treeop{+}{25}
    \rrtree{b}{b}{a}
    \\[-20pt]
    & \; \tiny
    \treeop{+}{25}
    \rrtree{b}{a}{b}
    \treeop{+}{25}
    \rrtree{b}{b}{b}
    \normalsize \treeop{= 0}{25}
    \\[-10pt]
    & \ingamma{\gamma_{2,3,1}\left(\right.}{4}
                                \hspace{-5pt}\treeingamma{\\[-6pt]b}\hspace{-10pt}
                                \ingamma{, \; (Q8') \; , \; \vert \;\left.\right)}{4}
                                \hspace{-20pt} \treeop{\quad - \quad}{4} \hspace{-20pt}
                                \ingamma{\gamma_{3,2,1,1}\left(\right.\; (Q7') \; ,}{4}
                                \hspace{-5pt}\treeingamma{\\[-6pt]a}\hspace{-5pt}
                                \ingamma{, \; \vert \; , \; \vert \;\left.\right)}{4}
    \\[-25pt]
    & \; \treeop{=}{25} \tiny
    \lltree{b}{b}{a}
    \treeop{+}{25}
    \lltree{b}{b}{d}_{(Q7')}
    \treeop{-}{25}
    \lrtree{b}{a}{b}_{(Q8')}
    \treeop{-}{25}
    \lltree{b}{b}{a}
    \treeop{-}{25}
    \lltree{b}{c}{a}_{(Q5')}
    \\[-20pt]
    & \; \tiny
    \treeop{+}{25}
    \corolla{b}{a}{a}_{(Q8')}
    \treeop{+}{25}
    \corolla{b}{a}{b}_{(Q8')}
    \\[-20pt]
    & \; \treeop{\equiv}{25} \tiny
    \treeop{-}{25}
    \lltree{b}{c}{d}_{(Q2')}
    \treeop{+}{25}
    \corolla{b}{d}{a}
    \treeop{+}{25}
    \corolla{b}{d}{b}
    \treeop{+}{25}
    \lrtree{b}{d}{b}
    \treeop{-}{25}
    \rltree{a}{b}{b}_{(Q7')}
    \\[-20pt]
    & \; \tiny
    \treeop{-}{25}
    \lrtree{b}{a}{c}_{(Q8')}
    \treeop{-}{25}
    \corolla{b}{d}{a}
    \treeop{+}{25}
    \rrtree{a}{b}{a}
    \treeop{-}{25}
    \corolla{b}{d}{b}
    \treeop{+}{25}
    \rrtree{a}{b}{b}
    \\[-20pt]
    & \; \treeop{\equiv}{25} \tiny
    \treeop{-}{25}
    \lrtree{b}{d}{b}
    \treeop{-}{25}
    \lrtree{b}{d}{c}
    \treeop{+}{25}
    \lrtree{b}{d}{b}
    \treeop{+}{25}
    \rltree{a}{b}{c}
    \treeop{-}{25}
    \rrtree{a}{b}{a}
    \treeop{-}{25}
    \rrtree{a}{b}{b}
    \\[-20pt]
    & \; \tiny
    \treeop{+}{25}
    \lrtree{b}{d}{c}
    \treeop{-}{25}
    \rltree{a}{b}{c}
    \treeop{+}{25}
    \rrtree{a}{b}{a}
    \treeop{+}{25}
    \rrtree{a}{b}{b}
    \normalsize \treeop{= 0}{25}
    \\[-10pt]
    & \ingamma{\gamma_{2,3,1}\left(\right.}{4}
                                \hspace{-5pt}\treeingamma{\\[-6pt]b}\hspace{-10pt}
                                \ingamma{, \; (Q9') \; , \; \vert \;\left.\right)}{4}
                                \hspace{-20pt} \treeop{\quad - \quad}{4} \hspace{-20pt}
                                \ingamma{\gamma_{3,2,1,1}\left(\right.\; (Q8') \; ,}{4}
                                \hspace{-5pt}\treeingamma{\\[-6pt]a}\hspace{-5pt}
                                \ingamma{, \; \vert \; , \; \vert \;\left.\right)}{4}
    \\[-25pt]
    & \; \treeop{=}{25} \tiny
    \lltree{b}{a}{a}
    \treeop{+}{25}
    \lltree{b}{a}{d}_{(Q8')}
    \treeop{+}{25}
    \lltree{b}{a}{b}_{(Q8')}
    \treeop{+}{25}
    \lltree{b}{a}{c}_{(Q8')}
    \treeop{-}{25}
    \lrtree{b}{a}{a}_{(Q8')}
    \\[-20pt]
    & \; \tiny
    \treeop{-}{25}
    \lltree{b}{a}{a}
    \treeop{-}{25}
    \lltree{b}{d}{a}_{(Q6')}
    \treeop{+}{25}
    \corolla{a}{a}{b}_{(Q9')}
    \\[-20pt]
    & \; \treeop{\equiv}{25} \tiny
    \treeop{-}{25}
    \lltree{b}{d}{d}_{(Q3')}
    \treeop{+}{25}
    \corolla{a}{d}{b}
    \treeop{-}{25}
    \lltree{b}{d}{b}
    \treeop{+}{25}
    \corolla{a}{b}{b}
    \treeop{-}{25}
    \lltree{b}{d}{c}
    \\[-20pt]
    & \; \tiny
    \treeop{+}{25}
    \corolla{a}{c}{b}
    \treeop{+}{25}
    \lrtree{b}{d}{a}
    \treeop{-}{25}
    \rltree{a}{b}{a}_{(Q8')}
    \treeop{+}{25}
    \lltree{b}{d}{b}
    \treeop{-}{25}
    \lrtree{b}{a}{d}_{(Q8')}
    \\[-20pt]
    & \; \tiny
    \treeop{-}{25}
    \corolla{a}{d}{b}
    \treeop{-}{25}
    \corolla{a}{b}{b}
    \treeop{-}{25}
    \corolla{a}{c}{b}
    \treeop{+}{25}
    \rrtree{a}{a}{b}
    \\[-20pt]
    & \; \treeop{\equiv}{25} \tiny
    \lltree{b}{d}{c}
    \treeop{-}{25}
    \lrtree{b}{d}{a}
    \treeop{-}{25}
    \lrtree{b}{d}{d}
    \treeop{-}{25}
    \lltree{b}{d}{c}
    \treeop{+}{25}
    \lrtree{b}{d}{a}
    \treeop{+}{25}
    \rltree{a}{b}{d}
    \\[-20pt]
    & \; \tiny
    \treeop{-}{25}
    \rrtree{a}{a}{b}
    \treeop{+}{25}
    \lrtree{b}{d}{d}
    \treeop{-}{25}
    \rltree{a}{b}{d}
    \treeop{+}{25}
    \rrtree{a}{a}{b}
    \normalsize \treeop{= 0}{25}
    \\[-10pt]
    & \ingamma{\gamma_{2,3,1}\left(\right.}{4}
                                \hspace{-5pt}\treeingamma{\\[-6pt]a}\hspace{-10pt}
                                \ingamma{, \; (Q9') \; , \; \vert \;\left.\right)}{4}
                                \hspace{-20pt} \treeop{\quad - \quad}{4} \hspace{-20pt}
                                \ingamma{\gamma_{3,2,1,1}\left(\right.\; (Q9') \; ,}{4}
                                \hspace{-5pt}\treeingamma{\\[-6pt]a}\hspace{-5pt}
                                \ingamma{, \; \vert \; , \; \vert \;\left.\right)}{4}
    \\[-25pt]
    & \; \treeop{=}{25} \tiny
    \lltree{a}{a}{a}
    \treeop{+}{25}
    \lltree{a}{a}{d}_{(Q9')}
    \treeop{+}{25}
    \lltree{a}{a}{b}_{(Q9')}
    \treeop{+}{25}
    \lltree{a}{a}{c}_{(Q9')}
    \treeop{-}{25}
    \lrtree{a}{a}{a}_{(Q9')}
    \\[-20pt]
    & \; \tiny
    \treeop{-}{25}
    \lltree{a}{a}{a}
    \treeop{-}{25}
    \lltree{a}{d}{a}_{(Q6')}
    \treeop{-}{25}
    \lltree{a}{b}{a}_{(Q8')}
    \treeop{-}{25}
    \lltree{a}{c}{a}_{(Q5')}
    \treeop{+}{25}
    \corolla{a}{a}{a}_{(Q9')}
    \\[-20pt]
    & \; \treeop{\equiv}{25} \tiny
    \treeop{-}{25}
    \lltree{a}{d}{d}_{(Q3')}
    \treeop{-}{25}
    \lltree{a}{b}{d}
    \treeop{-}{25}
    \lltree{a}{c}{d}_{(Q2')}
    \treeop{+}{25}
    \corolla{a}{d}{a}
    \treeop{-}{25}
    \lltree{a}{d}{b}
    \\[-20pt]
    & \; \tiny
    \treeop{-}{25}
    \lltree{a}{b}{b}_{(Q7')}
    \treeop{-}{25}
    \lltree{a}{c}{b}_{(Q4')}
    \treeop{+}{25}
    \corolla{a}{b}{a}
    \treeop{-}{25}
    \lltree{a}{d}{c}
    \treeop{-}{25}
    \lltree{a}{b}{c}
    \\[-20pt]
    & \; \tiny
    \treeop{-}{25}
    \lltree{a}{c}{c}_{(Q1')}
    \treeop{+}{25}
    \corolla{a}{c}{a}
    \treeop{+}{25}
    \lrtree{a}{d}{a}
    \treeop{+}{25}
    \lrtree{a}{b}{a}
    \treeop{+}{25}
    \lrtree{a}{c}{a}
    \\[-20pt]
    & \; \tiny
    \treeop{-}{25}
    \rltree{a}{a}{a}_{(Q9')}
    \treeop{+}{25}
    \lltree{a}{d}{b}
    \treeop{-}{25}
    \lrtree{a}{a}{d}_{(Q9')}
    \treeop{+}{25}
    \lltree{a}{b}{d}
    \treeop{-}{25}
    \lrtree{a}{a}{b}_{(Q9')}
    \\[-20pt]
    & \; \tiny
    \treeop{-}{25}
    \lrtree{a}{a}{c}_{(Q9')}
    \treeop{-}{25}
    \corolla{a}{d}{a}
    \treeop{-}{25}
    \corolla{a}{b}{a}
    \treeop{-}{25}
    \corolla{a}{c}{a}
    \treeop{+}{25}
    \rrtree{a}{a}{a}
    \\[-20pt]
    & \; \treeop{\equiv}{25} \tiny
    \lltree{a}{d}{c}
    \treeop{-}{25}
    \lrtree{a}{d}{a}
    \treeop{-}{25}
    \lrtree{a}{d}{d}
    \treeop{-}{25}
    \lrtree{a}{d}{b}
    \treeop{-}{25}
    \lrtree{a}{d}{c}
    \treeop{+}{25}
    \lltree{a}{b}{c}
    \\[-20pt]
    & \; \tiny
    \treeop{-}{25}
    \lrtree{a}{b}{a}
    \treeop{-}{25}
    \lrtree{a}{b}{b}
    \treeop{-}{25}
    \lrtree{a}{b}{d}
    \treeop{-}{25}
    \lrtree{a}{b}{c}
    \treeop{-}{25}
    \lltree{a}{d}{c}
    \treeop{-}{25}
    \lltree{a}{b}{c}
    \\[-20pt]
    & \; \tiny
    \treeop{-}{25}
    \lrtree{a}{c}{a}
    \treeop{-}{25}
    \lrtree{a}{c}{d}
    \treeop{-}{25}
    \lrtree{a}{c}{b}
    \treeop{-}{25}
    \lrtree{a}{c}{c}
    \treeop{+}{25}
    \lrtree{a}{d}{a}
    \treeop{+}{25}
    \lrtree{a}{b}{a}
    \\[-20pt]
    & \; \tiny
    \treeop{+}{25}
    \lrtree{a}{c}{a}
    \treeop{+}{25}
    \rltree{a}{a}{d}
    \treeop{+}{25}
    \rltree{a}{a}{b}
    \treeop{+}{25}
    \rltree{a}{a}{c}
    \treeop{-}{25}
    \rrtree{a}{a}{a}
    \treeop{+}{25}
    \lrtree{a}{d}{d}
    \\[-20pt]
    & \; \tiny
    \treeop{+}{25}
    \lrtree{a}{b}{d}
    \treeop{+}{25}
    \lrtree{a}{c}{d}
    \treeop{-}{25}
    \rltree{a}{a}{d}
    \treeop{+}{25}
    \lrtree{a}{d}{b}
    \treeop{+}{25}
    \lrtree{a}{b}{b}
    \treeop{+}{25}
    \lrtree{a}{c}{b}
    \\[-20pt]
    & \; \tiny
    \treeop{-}{25}
    \rltree{a}{a}{b}
    \treeop{+}{25}
    \lrtree{a}{d}{c}
    \treeop{+}{25}
    \lrtree{a}{b}{c}
    \treeop{+}{25}
    \lrtree{a}{c}{c}
    \treeop{-}{25}
    \rltree{a}{a}{c}
    \treeop{+}{25}
    \rrtree{a}{a}{a}
    \normalsize \treeop{= 0}{25}
    \\[-15pt]
\end{align*}
Since all compositions reduce to zero, we can state the following result.

\begin{theorem}
    The set of defining identities for quadri-algebras is a self-reduced Gr\"obner-Shirshov basis for the non-symmetric
    quadri-algebra operad with respect to the path-lexicographical ordering of monomials and with orders of operations
    $ c \; < \; b \; < \; d \; < \; a $ or $ c \; < \; d \; < \; b \; < \; a $.
\end{theorem}

Since the quadri-algebra operad is quadratic and its Gr\"obner-Shirshov basis is also quadratic,
from \cite[Prop.~3.10]{Hoffbeck10} and \cite[Corollary~3]{DotsenkoKhoroshkin10}, we immediately get
that the non-symmetric quadri-algebra operad is Koszul. Thus we found an alternative proof of this fact.

We can describe the normal tree monomials with respect to this Gr\"obner-Shirshov basis.
For normal tree monomials, any growth to the right is allowed, but the following growths to the left are not permitted \\[-15pt]
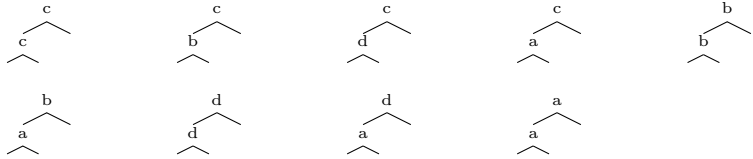
\begin{align*}
    {\tiny \ltree{c}{c} } & & {\tiny \ltree{c}{b} } & & {\tiny \ltree{c}{d} } & & {\tiny \ltree{c}{a} } & & {\tiny \ltree{b}{b} } \\
    {\tiny \ltree{b}{a} } & & {\tiny \ltree{d}{d} } & & {\tiny \ltree{d}{a} } & & {\tiny \ltree{a}{a} } & & \\[-15pt]
\end{align*}

The dimension of the multilinear subspace of degree $n$ of the free quadri-algebra in one variable
was conjectured by Aguiar and Loday \cite{AguiarLoday04} and proved by Vallette \cite{Vallette08} to be
\[
    d_n = \frac{1}{n} \sum_{j=n}^{2n-1} \binom{3n}{n+1+j} \binom{j-1}{j-n},
\]
which is the number of non-crossing connected graphs on $n+1$ vertices.


\section{Conclusions and open problems}

In this paper we described how to compute Gr\"obner-Shirshov bases for the non-symmetric operad of quadri-algebras and
found a quadratic one, which implies that this operad is Koszul. As a result, we obtained Gr\"obner-Shirshov bases for
the free quadri-algebra.
We also improved the result by Chen and Wang \cite{ChenWang10} by giving a shorter Gr\"obner-Shirshov basis for the free
dendriform algebra and a simpler proof using the operadic framework.

The defining identities for quadri-algebras are given in terms of the operations $\prec$, $\succ$, $\wedge$, $\vee$ and $*$.
We can then consider different (but equivalent) sets of operations and defining relations and compute the
corresponding Gr\"obner-Shirshov bases associated to the possible orders of operations. This means we change the $S_2$-module
$\mathcal{P}(2)$ of the operad and the generators of the operadic ideal $(R)$ of relators.

In \cite{BremnerMadariaga13}, special identities for the pre-Jordan algebras are found.
Pre-Lie and pre-Jordan algebras can be obtained from bilinear operations on dendriform algebras
analogously as Lie and Jordan algebras can be obtained from bilinear operations on associative algebras.
This process of obtaining nonassociative structures from associative structures can be extended to obtain
L-dendriform algebras \cite{BaiLiuNi10} and J-dendriform algebras \cite{BaiHou12} from quadri-algebras.
Having a Gr\"obner-Shirshov basis for quadri-algebras allows to extend the results of \cite{BremnerMadariaga13}
to the quadri-algebra setting.

If we iterate this splitting procedure we obtain a series of structures with $2^n$ binary nonassociative operations,
namely Loday algebras \cite{Vallette08} or ABQR algebras \cite{EbrahimiFardGuo05}; see also \cite{BaiBellierGuoNi12}.
It would be interesting to have Gr\"obner-Shirshov bases for these structures,
as well as for the algebras defined by the dual operads \cite{Vallette08}.


\section*{Acknowledgements}

The author thanks Prof. M.R.~Bremner for suggesting this problem to her, for his scientific generosity, 
his constant help and support and his useful comments and suggestions.
She also thanks the CIMAT in Guanajuato for holding the CIMPA research school
``Associative and Nonassociative Algebras and Dialgebras: Theory and Algorithms''
and Prof. Vladimir Dotsenko for his course on Gr\"obner bases for operads,
main source of inspiration for considering this problem, and his helpful remarks.


\end{document}